\documentclass[12pt]{article}

\usepackage{latexsym}
\usepackage{amssymb}
\usepackage{graphicx}
\usepackage{color}
\usepackage{amsmath}
\usepackage{mathrsfs} 

\newtheorem{theorem}{Theorem}

\newtheorem{corollary}[theorem]{Corollary}

\newtheorem{lemma}[theorem]{Lemma}

\newtheorem{remark}[theorem]{Remark}

\begin{document}

\title{Growth of subsolutions to fully nonlinear equations in halfspaces}

\author{
NIKLAS L.P. LUNDSTR\"OM
\\\\
\it \small Department of Mathematics and Mathematical Statistics,\\
\it \small Ume{\aa} University, SE-90187 Ume{\aa}, Sweden\/{\rm ;}\\
\it \small niklas.lundstrom@umu.se
}

\maketitle


\begin{abstract}

\noindent
We characterize lower growth estimates for
subsolutions in halfspaces of fully nonlinear partial differential equations on the form
$$
F(x,u,Du,D^2u) = 0
$$
in terms of solutions to ordinary differential equations
built solely upon a growth assumption on $F$.
Using this characterization we derive several sharp
Phragmen--Lindel\"of-type theorems for certain classes of well known PDEs.
The equation need not be uniformly elliptic nor homogeneous and
we obtain results both in case the subsolution is bounded or unbounded.
Among our results we retrieve classical estimates in the halfspace for $p$-subharmonic functions and extend those to more general equations;
we prove sharp growth estimates, in terms of $k$
and the asymptotic behaviour of
$\int_{0}^{R} C(s) ds$, for subsolutions
of equations allowing for sublinear growth in the gradient of the form
$C(|x|)|Du|^k$ with $k\geq 1$;
we establish a Phragmen--Lindel\"of theorem for weak subsolutions of
the variable exponent $p$-Laplace equation in halfspaces,
$1 < p(x) < \infty$, $p(x) \in C^1$,
of which we conclude sharpness by finding the ``slowest growing" $p(x)$-harmonic function together with its corresponding family of $p(x)$-exponents.
The paper ends with a discussion of our results from the point of view of a spatially dependent diffusion problem.\\

%
%
%
%
%
%
%
%
%
%
%
%
%
%
%
%
%
%
%
%
\end{abstract}


\section{Introduction}

\setcounter{theorem}{0}
\setcounter{equation}{0}

%
%
%
%

\noindent
We consider fully nonlinear nonhomogeneous elliptic partial differential equations in nondivergence form,
\begin{align}\label{eq:main}
F(x,u,Du,D^2u) = 0,\tag{$\star$}
\end{align}
in halfspaces in $\mathbb{R}^n$ where $n \geq 1$.
Here, $Du$ is the gradient, $D^2u$ the hessian,
$F:\mathbb{R}^n \times \mathbb{R} \times \mathbb{R}^n  \times \mathbb{S}^n \rightarrow \mathbb{R}$ and $\mathbb{S}^n$ is the set of symmetric $n \times n$ matrices equipped with the positive semi-definite ordering;
for $X, Y \in \mathbb{S}^n$, we write $X \leq Y$ if $\langle (X - Y) \xi, \xi \rangle \leq 0$ for all $\xi \in \mathbb{R}^n$.
Without loss of generality we fix the halfspace to $\mathbb{R}^n_+ := \{x \in \mathbb{R}^n : x_n > 0\}$ and assume the following:
\begin{itemize}
\item[] Degenerate ellipticity holds, i.e. $F(x,u,p,X) \geq  F(x,v,p,Y)$
whenever $u \geq v$, $X \leq Y$, as well as the growth condition
\begin{align}\label{eq:ass_drift_super}
-F(x,0,p,X) &\leq \Phi(|x|,|p|) + \Lambda(x_n) \text{Tr}(X^+) - \lambda(x_n) \text{Tr}(X^-)
\tag{$\star\star$}
\end{align}
whenever $x, p \in \mathbb{R}^n, X \in \mathbb{S}^n$, $X = X^+ - X^-$, $X^+ \geq 0$, $X^- \geq 0$ and $X^+  X^- = 0$.
Here,
$\Phi : [0,\infty) \times [0,\infty) \to (-\infty,\infty)$ is continuous,
nonincreasing in its first argument
 and
$\lambda, \Lambda : [0,\infty) \to (0,\infty)$ are functions such that
${\lambda}$ is nonincreasing and ${\Lambda}$ is nondecreasing.
\end{itemize}
Concerning $\Phi$ we will also need the following assumption:
\begin{itemize}
  \item[]
Either $\Phi$ is nonnegative and it holds that for all $\epsilon, t > 0$
(interpreting $1/0 = \infty$)
\begin{align}\label{eq:ass_phi_near}
\int_{0}^{\epsilon} \frac{ds}{\Phi(t, s)} = \infty
\tag{$\star\star\star$}
\end{align}
or $\Phi$ is nonpositive and $-\Phi(t,s)$ satisfies \eqref{eq:ass_phi_near} for all $\epsilon, t > 0$.
\end{itemize}

Under assumptions \eqref{eq:ass_drift_super}--\eqref{eq:ass_phi_near} we characterize the growth of
viscosity subsolutions of \eqref{eq:main} in halfspaces in terms of solutions to ODEs
(Theorem \ref{th:Phragmen--Lindelof}) 
 which are built solely upon functions $\Phi, \lambda$ and $\Lambda$ in \eqref{eq:ass_drift_super}.
Using this characterization we are able to derive sharp growth estimates of Phragmen--Lindel\"of type once the solutions to the ODEs are sufficiently understood.
Indeed, to apply Theorem \ref{th:Phragmen--Lindelof} one needs to
{\bf (1)} find functions $\Phi, \lambda$ and $\Lambda$ to ensure assumptions
\eqref{eq:ass_drift_super} and \eqref{eq:ass_phi_near},
{\bf (2)} solve the corresponding ODEs given in \eqref{eq:ODE} and
{\bf (3)} find the limit in Theorem \ref{th:Phragmen--Lindelof}.
An estimate is obtained if this limit is positive.
Theorem \ref{th:Phragmen--Lindelof} applies both in case the subsolution is bounded or unbounded,
and it can be used to find such border.

In Section \ref{sec:applic} we apply Theorem \ref{th:Phragmen--Lindelof} to derive sharp estimates for subsolutions to some well known PDEs of which the corresponding ODEs can be solved explicitly.
For example, we retrieve the classical Phragmen--Lindel\"of theorem in the halfspace for $p$-subharmonic functions by Lindqvist \cite{L85} and show in addition that it holds also for equations of $p$-Laplace type with lower order terms and vanishing ellipticity.
We obtain sharp lower estimates of the growth, in terms of $k \geq 1$
and the asymptotic behaviour of
$\int_{0}^{R} C(s)\lambda^{-1}(s) ds$, for subsolutions of equations with sublinear growth in the gradient such as
\begin{align*}
-P^-_{\lambda, \Lambda}(D^2u) + C(|x|)|Du|^k = 0
\end{align*}
in which $P^-_{\lambda, \Lambda}$ is a Pucci operator (definition recalled below) and $C(t)$ is a nonincreasing function.
These results reveal e.g. the border determining if a subsolution must grow to infinity or not in terms of $C(t)$ and $k$,
see Corollary \ref{cor:sublinear1} and estimate \eqref{eq:conclusion-sublinear}.
Moreover,
Theorem \ref{th:Phragmen--Lindelof} applies to nonhomogeneous PDEs including the variable exponent $p$-Laplace equation
$$
\nabla\cdot\left(|Du|^{p(x)-2} Du \right) = 0
$$
and we prove a sharp Phragmen--Lindel\"of theorem for weak subsolutions of this equation whenever $1 < p(x) < \infty$ is $C^1$ regular (Theorem \ref{thm:p(x)}).
It turns out that the growth estimate heavily depends on whether the subsolution
ever exceeds $x_n$ (distance to boundary) or not.
We conclude sharpness by finding the ``slowest growing" $p(x)$-harmonic function in the halfspace, for a given ellipticity bound, together with its corresponding family of $p(x)$ exponents (Remark \ref{cor:tjohej}). In the geometric setting of halfspaces,
this theorem sharpens some results of Adamowicz \cite{A14}.

The proof of Theorem \ref{th:Phragmen--Lindelof} relies on comparison with certain classical supersolutions of \eqref{eq:main}, which we construct in Lemma \ref{le:barrier_super} using solutions to the aforementioned ODEs.
We stress generality by pointing out that with the validity of Theorem \ref{th:Phragmen--Lindelof} at hand,
growth estimates for subsolutions to certain PDEs not considered in Section \ref{sec:applic} can be proved mainly by estimating solutions of first order ODEs and limits.

We end the paper by discussing the problem under investigation from the point of view of a diffusion problem.
Indeed, in Section \ref{sec:diffusion} we briefly explain,
through the application of spatially dependent diffusion,
why parts of our results presented in Theorem \ref{thm:p(x)} should hold.



We remark that our main results allow for ellipticity to blow up at infinity as $\lambda(x_n)$ may vanish and $\Lambda(x_n)$ may explode at infinity.
Moreover, the Osgood-type condition in \eqref{eq:ass_phi_near} is necessary to ensure that
subsolutions must continue to grow.
Indeed, for the strong maximum principle, see Julin \cite{J13}, Lundstr\"om--Olofsson--Toivanen \cite{LOT20} and the remarks below Theorem \ref{th:Phragmen--Lindelof}.
Furthermore, assumption \eqref{eq:ass_drift_super} can be written, with $\lambda = \lambda(x_n)$ and $\Lambda = \Lambda(x_n)$,
\begin{align*}
-F(x,0,p,X) \leq \Phi(|x|,|p|) - \mathcal{P}^-_{\lambda,\Lambda}(X) \quad \text{whenever} \quad x, p \in \mathbb{R}^n, X \in \mathbb{S}^n,
\end{align*}
where $\mathcal{P}^-_{\lambda,\Lambda}(X) = -\Lambda \text{Tr}(X^+) + \lambda \text{Tr}(X^-)$ is the Pucci maximal operator, $X = X^+ - X^-$ with $X^+ \geq 0$, $X^- \geq 0$ and $X^+  X^- = 0$.
In particular, if $X \in \mathbb{S}^n$ has eigenvalues $e_1,e_2, \dots ,e_n$
the Pucci extremal operators $\mathcal{P}^+_{\lambda,\Lambda}$ and $\mathcal{P}^-_{\lambda,\Lambda}$ with ellipticity $0 < \lambda \leq \Lambda$ are defined by
$$
\mathcal{P}^+_{\lambda,\Lambda}(X) := - \lambda \sum_{e_i \geq 0} e_i - \Lambda\sum_{e_i < 0} e_i \quad \text{and} \quad \mathcal{P}^-_{\lambda,\Lambda}(X) := - \Lambda \sum_{e_i \geq 0} e_i - \lambda \sum_{e_i < 0} e_i.
$$
For properties of the Pucci operators see e.g. Caffarelli--Cabre \cite{CC95} or Capuzzo-Dolcetta--Vitolo \cite{CDV07}.
We remark also that
the above assumption \eqref{eq:ass_drift_super} is implied by the standard ellipticity assumption 
\begin{align}\label{eq:ass_ellipt}
\lambda \text{Tr}(Y) \leq F(x,u,p,X) - F(x,u,p,X + Y) \leq \Lambda \text{Tr}(Y), 
\end{align}
whenever $Y$ is positive semi-definite, together with
\begin{align}\label{eq:ass_drift_super_old}
-F(x,0,p,0) \leq \Phi(|x|, |p|) \quad \text{whenever} \quad x, p \in \mathbb{R}^n.
\end{align}
Observe also that \eqref{eq:ass_drift_super} allows for nonlinear degenerate
elliptic operators which do not satisfy \eqref{eq:ass_ellipt}.
For example operators of the form
$$
F(X) = - \Lambda \left( \sum_{i=1}^{n} \Gamma(\mu_i^+)\right) + \lambda \left( \sum_{i=1}^{n} \Psi(\mu_i^-) \right)
$$
where $\mu_i, i = 1,\dots,n$, are the eigenvalues of the matrix $X \in \mathbb{S}^n$ and $\Gamma, \Psi : [0,\infty) \to [0,\infty)$
are continuous and nondecreasing functions such that $\Gamma(s) \leq s \leq \Psi(s)$, see Capuzzo-Dolcetta--Vitolo \cite{CDV07}.


The Phragm\'en-Lindel\"of principle and results of Phragm\'en-Lindel\"of type,
which has connections to elasticity theory
(Horgan \cite{Horgan}, Quintanilla \cite{Q93}, Leseduarte--Carme--Quintanilla \cite{Qnew}),
have been frequently studied during the last century.
To mention a few papers (without giving a complete summary),
Ahlfors \cite{A37} extended results from Phragm\'en--Lindel\"of \cite{PL08} to the upper half space of $\mathbb{R}^n$,
Gilbarg \cite{G52}, Serrin \cite{S54} and Herzog \cite{H64} considered more general elliptic equations of second order.
Miller \cite{M71} considered uniformly elliptic operators in nondivergence form and
unbounded domains contained in cones.
Kurta \cite{K93} and Jin--Lancaster \cite{JL03} estimated growth of bounded solutions of quasilinear equations, the later used solutions to boundary value problems, while
Vitolo \cite{V04} considered elliptic equations in sectors.
Capuzzo-Dolcetta--Vitolo \cite{CDV07} and Armstrong--Sirakov--Smart \cite{ASS12} considered fully nonlinear equations, the later in certain Lipschitz domains, and
Koike--Nakagawa \cite{KN09} established Phragm\'en-Lindel\"of theorems for subsolutions of fully nonlinear elliptic PDEs with unbounded coefficients and inhomogeneous terms.
Adamowicz \cite{A14} studied subsolutions of the variable exponent $p$-Laplace equation,
while Bhattacharya \cite{Bhatt05} and Granlund--Marola \cite{GM14} considered infinity-harmonic functions.
Lindqvist \cite{L85} established Phragm\'en-Lindel\"of's theorem for $n$-subharmonic functions when the boundary is an $m$-dimensional hyperplane in $\mathbb{R}^n$,
$0 \leq m \leq n-1$, which was extended to $p$-subharmonic functions,
$n-m < p \leq \infty$, in Lundstr\"om \cite{L16}.
We also mention that recently, Braga--Moreira \cite{BM20} showed that nonnegative solutions to a generalized $p$-Laplace equation in the upper halfplane, vanishing on $\{x_n = 0\}$, is $u(x)= x_n$ (modulo normalization) and
Lundstr\"om--Singh \cite{LS21} proved a similar result for $p$-harmonic functions in planar sectors as well as a sharp Phragmen--Lindel\"of theorem.
Lundberg--Weitsman \cite{LW15} studied the growth of solutions to the minimal surface equation over domains containing a halfplane.
The spatial behavior of solutions of
the Laplace equation on a semi-infinite cylinder with dynamical nonlinear
boundary conditions was investigated in
Leseduarte--Carme--Quintanilla \cite{Qnew}.
Finally, we mention that recently, local estimates such as a sharp Harnack inequality (Julin \cite{J13}), Boundary Harnack inequalities (Avelin--Julin \cite{AJ17}) as well as
strong maximum and minimum principles (Lundstr\"om--Olofsson--Toivanen \cite{LOT20}) were established for fully nonlinear PDEs covered by the class of equations considered here.



\subsection*{Preliminaries}
\label{sec:prel}

For a point $x \in \mathbb{R}^n$ we use the notation $x = (x_1, x_2, \dots x_{n-1}, x_n) = (x', x_n)$.
By $\Omega$ we denote a domain, that is, an open connected set.
For a set $E \subset \mathbb{R}^n$ we let $\overline{E}$ denote the closure and $\partial E$ the boundary of $E$.
By $c$ we denote a positive constant not necessarily the same at each occurrence.
We write $A \precsim B$ if there exists $c$ such that $A \leq cB$.

A function $u : \Omega \to \mathbb{R}$ is a classical subsolution (supersolution) to \eqref{eq:main} in $\Omega$ if it is twice differentiable in $\Omega$ and satisfies $F(x,u,Du,D^2u) \leq 0$ ($F(x,u,Du,D^2u) \geq 0$).
If the inequality holds strict then $u$ is a strict classical subsolution (supersolution),
and if equality holds then it is a classical solution.

We choose to present our main results for viscosity subsolutions, of which we recall the definition below
in case $F:\mathbb{R}^n \times \mathbb{R} \times \mathbb{R}^n  \times \mathbb{S}^n \rightarrow \mathbb{R}$ and $F$ is a continuous function (which is not necessary for our results).

The following definition is from Crandall--Ishii--Lions \cite{CIL92}:
An upper semicontinuous (USC) function $u : \Omega \to \mathbb{R}$ is a viscosity subsolution 
if for any $\varphi \in C^2(\Omega)$ and any $x_0 \in \Omega$
such that $u - \varphi$ has a local maximum at $x_0$
it holds that
\begin{align*}
F(x_0,u(x_0),D\varphi(x_0),D^2\varphi(x_0)) \leq 0.
\end{align*}
A lower semicontinuous (LSC) function $u : \Omega \to \mathbb{R}$ is a viscosity supersolution 
if for any $\varphi \in C^2(\Omega)$ and any $x_0 \in \Omega$
such that $u - \varphi$ has a local minimum at $x_0$
it holds that
\begin{align*}
F(x_0,u(x_0),D\varphi(x_0),D^2\varphi(x_0)) \geq 0.
\end{align*}
A continuous function is a viscosity solution if it is both a viscosity sub- and a viscosity supersolution.

Let $u$ be a subsolution and $v$ a supersolution to \eqref{eq:main} and let $a$ and $b$ be constants.
As \eqref{eq:main} is not necessarily homogeneous,  $a + b u$ and $a + b v$ may fail as sub- and supersolutions.
However, degenerate ellipticity guarantees that $u - c$ is a subsolution,
and $u + c$ is a supersolution whenever $c \geq 0$.

We will not discuss the validity of a general comparison principle for viscosity solutions of \eqref{eq:main}
since we only need the possibility to compare viscosity subsolutions to  classical supersolutions which is possible.
Indeed, let $\Omega$ be a bounded domain,
$u$ a viscosity subsolution and $v$ a classical strict supersolution in $\Omega$,
$u \leq v$ on $\partial \Omega$ and that $u \geq v$ somewhere in $\Omega$.
By USC the function $u-v$ attains a maximum at some point $x_0 \in \Omega$.
Since $v \in C^2(\Omega)$, $u - v$ has a maximum at $x_0$ and $u$ is a viscosity subsolution it follows by definition of viscosity solutions that
\begin{align}\label{eq:johej-lemma1}
F(x_0, u(x_0), Dv(x_0), D^2 v(x_0)) \leq 0.
\end{align}
But since $v$ is a classical strict supersolution we have $F(x,v(x), Dv(x), D^2v(x)) > 0$ whenever $x \in \Omega$,
and as $u(x_0) \geq v(x_0)$ it follows from degenerate ellipticity that
%
$
F(x_0, u(x_0), Dv(x_0), D^2 v(x_0)) \geq F(x_0, v(x_0), Dv(x_0), D^2 v(x_0)) > 0.
$
%
This contradicts \eqref{eq:johej-lemma1} and hence we have proved the following simple lemma:

\begin{lemma}\label{le:comp-weak}
Let $\Omega$ be a bounded domain, $u \in USC(\overline{\Omega})$ a viscosity subsolution and $v \in LSC(\overline{\Omega})$
a viscosity supersolution of \eqref{eq:main} in $\Omega$ satisfying $u \leq v$ on $\partial \Omega$.
Assume degenerate ellipticity.
If either $u$ is a strict classical subsolution, or $v$ is a strict classical supersolution,
then $u < v$ in $\Omega$.
\end{lemma}

Neither the choice of viscosity solutions nor the assumption that $F$ is continuous are necessary for our results.
Any other definition of ``weak solutions" can be considered,
whenever more appropriate for the equation, as long as such weak subsolutions of \eqref{eq:main} can be compared to classical strict supersolutions of \eqref{eq:main}.
In particular, 
our proof relies on construction of a classical strict supersolution to \eqref{eq:main} and
comparison with this barrier function.
What is needed is the validity of following simple comparison result:

\begin{lemma}\label{le:comp-weak-some-sence}
Let $\Omega$ be a bounded domain, $u \in USC(\overline{\Omega})$ a subsolution (in some weak sense)
and $v$ a classical strict supersolution to \eqref{eq:main} in $\Omega$, continuous on $\overline{\Omega}$.
If $u \leq v$ on $\partial \Omega$ then $u \leq v$ in $\Omega$.
\end{lemma}


\section{Characterizing growth in terms of solutions to ordinary differential equations}
\label{sec:res}
\setcounter{theorem}{0}
\setcounter{equation}{0}

We will estimate the growth of subsolutions to \eqref{eq:main} in terms of solutions
$f : [0, \infty) \to \mathbb{R}$ to the following initial value problems,
originating from assumption \eqref{eq:ass_drift_super}:
If $\Phi(t,s) \geq 0$ for all $t,s \in \mathbb{R}_+$ then we will make use of solutions to
\begin{align}\label{eq:ODE}
\frac{df}{dt} = -  \frac{\Phi(t,f(t))}{\lambda(t)} - K(R) \frac{\Lambda(t)}{\lambda(t)} f(t),  \quad t \in (0,R) \quad \text{with} \quad f(0) = \nu,
\end{align}
%
%
%
where $\nu \geq 0$ and $R > 0$.
Through the paper, we will by $f_{\nu,R} = f_{\nu,R}(t)$ denote the solution to \eqref{eq:ODE} with $K(R) = \frac{n}{\gamma(R)}$, in which $\gamma(R)$ appears in the domain defined in \eqref{eq:domain_2} below.
Further, we denote by $f_{\nu} = f_{\nu}(t)$ the solution of \eqref{eq:ODE} when $K \equiv 0$.
If $\Phi \leq 0$ for all $t,s \in \mathbb{R}_+$ then we use instead solutions of \eqref{eq:ODE} but with $\lambda(t)$ replaced by $\Lambda(t)$
in the first term on the right hand side of \eqref{eq:ODE}. 
We allow ourselves to simplify notation according to $\lambda = \lambda(\cdot), \Lambda = \Lambda(\cdot), K = K(R)$ and $\gamma = \gamma(R)$ whenever appropriate.

Since $\Phi(t, 0) = 0$ and $\nu \in [0, \infty)$ the solutions $f_{\nu,R}$ and $f_{\nu}$ will be nonegative.
Moreover, if $\Phi$ satisfies the Osgood-type condition \eqref{eq:ass_phi_near}
then the solutions will, for any $\nu > 0$, remain positive.
(The only nonpositive solutions are the trivial solutions $f_{\nu,R}\equiv f_{\nu}\equiv 0$ starting at $\nu = 0$.)
This plays a role in our main results, as pointed out in the remarks made below Theorem \ref{th:Phragmen--Lindelof}.
In Figures \ref{fig:s^k} and \ref{fig:variablep} several solutions of \eqref{eq:ODE} are plotted for some choices of $\Phi$.

\begin{figure}[!hbt]
\begin{center}
\includegraphics[height = 8cm, width = 13 cm, viewport=10 150 375 365,clip]{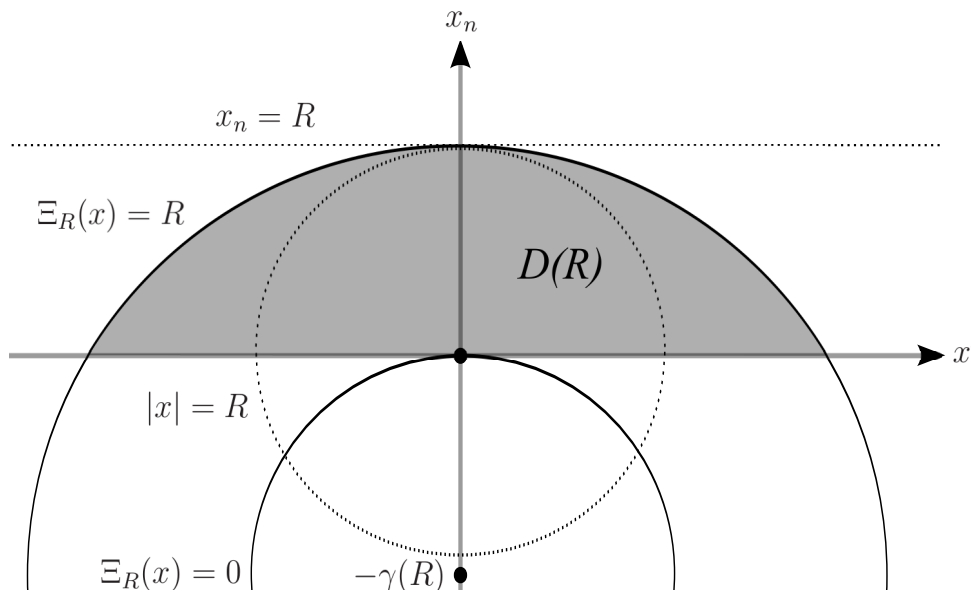}
\end{center}
\caption{Geometric definitions and constructions.}
\label{fig:geometry}
\end{figure}

To proceed we define, for a nondecreasing function $\gamma = \gamma(R) > 0$ and $n \geq 1$, the domain
\begin{align}\label{eq:domain_2}
D(R) :=  \left\{ x \in \mathbb{R}^n_+ :  \sum_{i=1}^{n-1} x_i^2 + (x_n + \gamma)^2 < (R + \gamma)^2  \right\},
\end{align}
%
see Figure \ref{fig:geometry}.
Finally, for a subsolution $u$ and for $R > 0$ we define
\begin{align*}
M(R) = \sup_{\partial D(R)} u 
\end{align*}
and
\begin{align*}
M'(R) = \liminf_{h \to 0^+} \frac{M(R) - M(R-h)}{h}.
\end{align*}

The following theorem characterizes a sharp lower growth estimate of subsolutions to
\eqref{eq:main} in terms of solutions $f_{\nu,R}$ and $f_{\nu}$ to the ODE \eqref{eq:ODE}:
\begin{theorem}
\label{th:Phragmen--Lindelof}
Suppose that \eqref{eq:ass_drift_super} and \eqref{eq:ass_phi_near} hold and let $u$ be a subsolution of \eqref{eq:main} in $\mathbb{R}^n_+$ satisfying
\begin{align*}
\limsup_{x \to y} \, u(x) \, \leq \,  0 \quad \textrm{for all} \quad y \in \partial \mathbb{R}^n_+. 
\end{align*}
Then either $u \leq 0$ in $\mathbb{R}^n_+$ or $M(R)$ is increasing and it holds that
\begin{align*}
\liminf_{R\to\infty} \, \frac{M'(R)}{f_{\nu}(R)} \, \geq \,\liminf_{R\to\infty} \, \frac{f_{\nu, R}(R)}{f_{\nu}(R)},
\end{align*}
where $\nu$ satisfies $u(\bar x) \geq \int_{0}^{\bar x_n} f_{\nu}(t)\, dt$ for some $\bar{x}$ on the $x_n$-axis.
\end{theorem}


\noindent
Using Theorem \ref{th:Phragmen--Lindelof} an ``explicit" growth estimate can thus be found
by estimating the limit $f_{\nu,R}(R)/f_{\nu}(R)$ as $R \to \infty$.
In Section \ref{sec:applic} we will consider certain PDEs for which we can solve the ODE \eqref{eq:ODE} explicitly -- calculate the limit -- and thereby prove Phragmen--Lindel\"of type theorems.
Let us note that if we can prove
$$
\liminf_{R\to\infty} \, \frac{M'(R)}{f_{\nu}(R)}  > 0
$$
then
$
M(R) - M(R_0) \geq c \int_{R_0}^{R} f_{\nu}(s) ds
$
whenever $R > R_0$ for some small $c$ and thus 
$$
\liminf_{R\to\infty} \, \frac{M(R)}{\int_{0}^{R} f_{\nu}(s) ds} > 0.
$$
Hence, if the integral
\begin{align*}
\int_{0}^{\infty} f_{\nu}(s) ds
\end{align*}
converges, then subsolutions may be bounded,
but if the integral diverges,
then subsolutions must grow to infinity and the conclusion of Theorem \ref{th:Phragmen--Lindelof} takes the form of classical Phrgmen--Lindel\"of theorems.

We remark that the assumption ``$\bar{x}$ lies on the $x_n$-axis" is only for notational simplicity;
we may translate coordinates otherwise.
Note also that Theorem \ref{th:Phragmen--Lindelof} holds whenever $\mathbb{R}^n_+$ is replaced (in the theorem and in \eqref{eq:domain_2}) with
$\Omega \subset\mathbb{R}^n_+$,
and that Theorem \ref{th:Phragmen--Lindelof} gives a growth estimate for any initial condition $\nu \geq 0$ in \eqref{eq:ODE}
as long as $u(\bar{x}) \geq \int_{0}^{\bar x_n} f_{\nu}(t)\, dt$ for some $\bar{x}$.
The best estimate corresponds to the largest $\nu$.
Moreover, it can be realized from the proof of Lemma \ref{le:barrier_super} that the assumption
``${\lambda}$ nonincreasing and ${\Lambda}$ nondecreasing" can be replaced by the slightly weaker assumption that
${\lambda}/{\Lambda}$ is nonincreasing and $\lambda$ is nonincreasing ($\Lambda$ is nondecreasing) when $\Phi\geq 0$ ($\Phi\leq 0$).
Finally, it will be clear from the proof that we also have
$
M(R) - M(R-h) \geq \int_{R}^{R+h} f_{\nu,R}(t) dt
$
for $R > \bar{x}_n$ and any $h > 0$.
We realize that $M(R)$ must increase as long as $f_{\nu,R} > 0$, which happens
whenever $\Phi$ satisfies the Osgood-type condition \eqref{eq:ass_phi_near}.
Otherwise, the strong maximum principle does not hold and a positive subsolution to \eqref{eq:main} may stop growing and attain an interior maximum, see Julin \cite{J13} and Lundstr\"om-Olofsson-Toivannen \cite[Remark 4.3]{LOT20} for a counterexample.


Concerning sharpness of Theorem \ref{th:Phragmen--Lindelof} we
consider the function
\begin{align}\label{eq:candidate}
u(x) = \int_{0}^{x_n} f_{\nu}(t) dt,
\end{align}
vanishing on $\partial\mathbb{R}_+^n$, depending only on $x_n$ with derivative $u_{x_n}(x) =  f_{\nu}(x_n)$.
In case $\Phi(t,s) \geq 0$ it holds that $u_{x_n x_n}(x) =  f'_{\nu}(x_n) = - \lambda^{-1}(x_n)  \Phi(x_n, f_{\nu}(x_n))$
and hence we obtain, e.g., that
\begin{align}\label{eq:candidatesPDE}
-F(x, Du, D^2u) :=  \lambda(x_n) \Delta u + \phi(f_{\nu}(x_n)) = 0,
\end{align}
for some function $\phi(s)$ satisfying \eqref{eq:ass_phi_near}.
In case $\Phi(t,s) \leq 0$ the same holds but with $\lambda$ replaced by $\Lambda$.
Thus, the function defined in \eqref{eq:candidate} is a classical solution to an equation of type \eqref{eq:main} satisfying \eqref{eq:ass_drift_super} and \eqref{eq:ass_phi_near}.
Moreover, it satisfies $M'(R) = f_{\nu}(R)$.
In conclusion, when the limit in Theorem \ref{th:Phragmen--Lindelof} is positive, the growth estimate cannot be improved, ignoring the shape of $D(R)$ and the value of the limit.

Concerning the shape of $D(R)$ we note the following.
If $n = 1$ then $D(R) = (0,R)$ independent of $\gamma$,
but if $n \geq 2$ then $\gamma = c R$ implies that the spherical segment $D(R)$ preserves its geometric proportions for all $R > 0$.
If $\gamma(R)/R$ is increasing then $D(R)$ expands faster in the $x'$-direction,
implying slightly weaker estimates since $\partial D(R)$, on which supremum is taken, becomes larger.
Observe that if the problem is considered in $\Omega \subset \mathbb{R}_+^n$ this might
be of minor importance, especially if e.g. $\Omega$ is bounded in $x'$-directions or contained in a cone with apex at the origin.
There is not much of a gain to take $\gamma(R)/R$ decreasing since $D(R)$ still expands at rate $R$ in $x'$-directions.

The proof of Theorem \ref{th:Phragmen--Lindelof} relies on comparison arguments
and the following construction of
a classical strict supersolution to \eqref{eq:main}.


\begin{lemma}
\label{le:barrier_super}
Suppose that \eqref{eq:ass_drift_super} and \eqref{eq:ass_phi_near} hold, let $R > 0$ and put
$$
\Xi_R(x) = \sqrt{ \sum_{i=1}^{n-1} x_i^2 + (x_n + \gamma)^2} - \gamma = |(x', x_n+\gamma)| - \gamma
$$
in which $\gamma = \gamma(R)$ is from \eqref{eq:domain_2}.
Then the function
\begin{align*}
V_R(x) =  \int_{0}^{\Xi_R(x)} f_{\nu,R}(t) dt
\end{align*}
is a strict classical supersolution to \eqref{eq:main} in $D(R)$. 
\end{lemma}


\noindent
{\bf Proof of Lemma \ref{le:barrier_super}.}
For notational simplicity we set
$\Xi = \Xi_R(x)$, $V = V_R(x)$ and $f(t) = f_{\nu,R}(t)$.
Differentiating yields
\begin{align*}
\frac{\partial V}{\partial x_i} = \frac{x_i}{|(x', x_n+\gamma)|} f\left(\Xi\right), \quad 1 \leq i \leq n-1, \quad
\frac{\partial V}{\partial x_n} = \frac{x_n + \gamma}{|(x', x_n+\gamma)|} f\left(\Xi\right).
\end{align*}
It follows that
\begin{align}\label{eq:grad-beg}
\vert D V \vert = f\left(\Xi\right).
\end{align}

\noindent
The second derivatives become
\begin{align*}
\frac{\partial^2 V}{\partial x_i^2}
= 
 \left( \frac{x_i}{|( x', x_n+\gamma)|} \right)^2 f'(\Xi)
+
  \left(\frac{1}{|(x', x_n+\gamma)|} - \frac{x_i^2}{|(x', x_n+\gamma)|^3} \right) f\left(\Xi\right),
%
\end{align*}
for $1 \leq i \leq n-1$, and
\begin{align*}
\frac{\partial^2 V}{\partial x_n^2} &=  \left(\frac{x_n + \gamma}{|(x', x_n+\gamma)|}\right)^2 f'(\Xi)
+
 \left(\frac{1}{|(x', x_n+\gamma)|} - \frac{(x_n + \gamma)^2}{|(x', x_n+\gamma)|^3} \right) f\left(\Xi\right),
\end{align*}
giving
\begin{align*}
\text{Tr}(D^2 V) 
=  f'(\Xi) + \frac{n-1}{|(x', x_n+\gamma)|}f\left(\Xi\right).
\end{align*}
By construction we have $f'(t) = - \frac{\Phi(t, f(t))}{\lambda(t)} - K \frac{\Lambda(t)}{\lambda(t)}f(t)$ and hence
\begin{align*}
\text{Tr}(D^2 V) 
= -\frac{\Phi\left(\Xi,f\left(\Xi\right)\right)}{\lambda(\Xi)} - K \frac{\Lambda(\Xi)}{\lambda(\Xi)} f(\Xi) +   \frac{n-1}{|(x', x_n+\gamma)|}f\left(\Xi\right).
\end{align*}
We assume from here on that $\Phi\geq 0$ and decompose $D^2 V = \left(D^2 V \right)^+ - \left(D^2 V \right)^-$ so that
%
\begin{align*}
\text{Tr}\left( \left(D^2 V \right)^+ \right) &=  \frac{n-1}{|(x', x_n+\gamma)|}f\left(\Xi\right) \quad \textrm{and} \\
\text{Tr}\left( \left(D^2 V \right)^- \right) &=  \frac{\Phi\left(\Xi,f\left(\Xi\right)\right)}{\lambda(\Xi)} + K \frac{\Lambda(\Xi)}{\lambda(\Xi)} f(\Xi)
\end{align*}
Utilizing the structure assumption \eqref{eq:ass_drift_super}, the fact that $V \geq 0$ and using \eqref{eq:grad-beg} give
\begin{align}\label{eq:super_proof_2}
F(x,V,DV,D^2V) \geq &F(x,0,DV,D^2V) \notag\\
\geq &- \Phi\left(|x|, f\left(\Xi\right)\right)
- \Lambda(x_n) \frac{n-1}{|(x', x_n+\gamma)|}f\left(\Xi\right)\notag\\
&+ \frac{\lambda(x_n)\Phi\left(\Xi,f\left(\Xi\right)\right)}{\lambda(\Xi)} + K \frac{\lambda(x_n)\Lambda(\Xi)}{\lambda(\Xi)} f(\Xi)   \notag\\
\geq  &- \Lambda(x_n)  \frac{n-1}{|(x', x_n+\gamma)|}f(\Xi)
+   K \frac{\lambda(x_n)\Lambda(\Xi)}{\lambda(\Xi)} f(\Xi)
\end{align}
since $\lambda(\Xi) \leq \lambda(x_n)$ and $\Phi\left(|x|, f\left(\Xi\right)\right) \leq \Phi\left(\Xi,f\left(\Xi\right)\right)$ holds.
This last statement follows since $|x| \geq \Xi_R(x) \geq x_n$ by geometry, see Figure \ref{fig:geometry}, 
and functions are nonincreasing by assumption.

To show that $V$ is a strict classical supersolution we need
$F(x,V,DV,D^2V) > 0$ and by \eqref{eq:super_proof_2} and the fact that $f(\Xi) > 0$ it suffices to ensure
$$
\frac{n-1}{|(x', x_n+\gamma)|} <  K \frac{\lambda(x_n)\Lambda(\Xi)}{\Lambda(x_n)\lambda(\Xi)}.
$$
Observing that
$\frac{\lambda(x_n)}{\Lambda(x_n)} \geq \frac{\lambda(\Xi)}{\Lambda(\Xi)}$ holds since $\lambda/\Lambda$ is nonincreasing,
it suffices to ensure
$$
\frac{n-1}{|(x', x_n+\gamma)|} <  K .
$$
We know that $\gamma \leq |(x', x_n+\gamma)|$ in $\mathbb{R}_+^n$  so it is enough to have
$$
\frac{n-1}{\gamma(R)} <  K
$$
which holds since we have $K = \frac{n}{\gamma(R)}$.

If $\Phi \leq 0$ then by construction $f'(t) = -\frac{\Phi\left(t,f\left(t\right)\right)}{\Lambda(t)} -  K \frac{\Lambda(t)}{\lambda(t)} f(t)$ and we obtain
\begin{align}\label{eq:tracetjoho}
\text{Tr}\left( \left(D^2 V \right)^+ \right) &= -\frac{\Phi\left(\Xi,f\left(\Xi\right)\right)}{\Lambda(\Xi)} + \frac{n-1}{|(x', x_n+\gamma)|} f(\Xi),\notag\\ 
\text{Tr}\left( \left(D^2 V \right)^- \right) &= K \frac{\Lambda(\Xi)}{\lambda(\Xi)} f(\Xi)
\end{align}
and thus instead of \eqref{eq:super_proof_2} we end up with
\begin{align}\label{eq:super_proof_33}
F(x,V,DV,D^2V) \geq &F(x,0,DV,D^2V)\notag\\
\geq &- \Phi\left(|x|, f\left(\Xi\right)\right)
 + \frac{\Lambda(x_n)}{\Lambda(\Xi)} \Phi\left(\Xi,f\left(\Xi\right)\right) \notag\\
&-\Lambda(x_n)\frac{n-1}{|(x', x_n+\gamma)|} f(\Xi) + \lambda(x_n) K \frac{\Lambda(\Xi)}{\lambda(\Xi)} f(\Xi) \notag\\
\geq &-\Lambda(x_n)\frac{n-1}{|(x', x_n+\gamma)|} f(\Xi) + \lambda(x_n) K \frac{\Lambda(\Xi)}{\lambda(\Xi)} f(\Xi).
\end{align}
Here, the last inequality holds since $x_n \leq \Xi_R(x) \leq |x|$,
$\Lambda$ is nondecreasing and $-\Phi \geq 0$ is nondecreasing in its first argument so that
$$
- \Phi\left(|x|, f\left(\Xi\right)\right) \geq \frac{\Lambda(x_n)}{\Lambda(\Xi)}\left( - \Phi\left(\Xi,f\left(\Xi\right)\right)\right).
$$
To ensure that $V$ is a strict supersolution we see from \eqref{eq:super_proof_33} that it remains to show that
$$
\frac{n-1}{|(x', x_n+\gamma)|} < K \frac{\lambda(x_n)}{\Lambda(x_n)} \frac{\Lambda(\Xi)}{\lambda(\Xi)}
$$
and we are thus back in the same situation as in the case $\Phi \geq 0$. 
%
%
The proof of Lemma \ref{le:barrier_super} is complete.
$\hfill \Box$ \\


\noindent
{\bf Proof of Theorem \ref{th:Phragmen--Lindelof}.}
Let $u$ be as in the statement of the theorem and denote with $\nu_0$ the initial condition in \eqref{eq:ODE} for which
we want to prove the growth estimate.
Let $R > 0$ and put $V := V_{R}(x)$,
where $V_{R}(x)$ is the strict supersolution in $D(R)$ guaranteed by Lemma \ref{le:barrier_super}.

If $V \geq M(R)$ on $\partial D(R)$ then, if $\Phi$ is nonnegative, it follows that $f_{\nu} \leq \nu$ and we obtain equality by decreasing $\nu$.
If $\Phi$ is nonpositive we note that $f'_{\nu,R}(t) \leq \widetilde\Phi(f_{\nu,R}(t))$ for some $\widetilde{\Phi}(s)\geq 0$ satisfying \eqref{eq:ass_phi_near}.
Thus
$$
\int_{\nu}^{f_{\nu,R}(t)}\frac{ds}{\widetilde{\Phi}(s)} \leq t
$$
which implies that $f_{\nu,R}(t)  \searrow 0$ as $\nu \searrow 0$ for all $t \in [0,R]$.
Therefore, we obtain equality by decreasing $\nu$ also in this case.
If $V < M(R)$ on $\partial D(R)$ then we increase $\nu$.
If this does not help,
(note that we may have $V \leq A$, for all $\nu$ and all $R$, for some $A > 0$),
then we lift the supersolution by adding a nonnegative constant.
Indeed, for $C \geq 0$ it follows from degenerate ellipticity that also $V + C$ is a strict supersolution.
We conclude that
$$
V + C \geq C \quad \text{on} \quad \{x_n = 0\}  \quad \text{and} \quad V + C = M(R)\quad \text{on} \quad \partial D(R) \cap \mathbb{R}^n_+.
$$
We clarify that if $C > 0$ then we have taken $\nu > \nu_0$.
It follows that
\begin{align*}
\limsup_{x \to z} u(x) \leq V(z) + C \quad \textrm{for all} \quad z \in \partial D(R)
\end{align*}
and the weak comparison principle in Lemma \ref{le:comp-weak}
implies that $u \leq V + C$ in $D(R)$.

We next conclude that $\nu \geq \nu_0$.
In particular, assume $\nu < \nu_0$.
By assumption and by the above we have
$V(\bar{x}) \geq  u(\bar{x}) \geq \int_{0}^{\bar x_n} f_{\nu_0}(t)\, dt$ for some $\bar{x} \in D(R) \cap \{\bar x' = 0\}$,
but on the other hand
\begin{align*}
V(\bar{x}) = \int_{0}^{\Xi_R(\bar x)} f_{\nu,R}\left(t\right) dt =
 \int_{0}^{\bar x_n} f_{\nu,R}\left( t\right) dt < \int_{0}^{\bar{x}_n} f_{\nu_0,R}\left(t\right) dt
 \leq \int_{0}^{\bar x_n} f_{\nu_0}(t)\, dt,
\end{align*}
where the last inequality follows since $f_{\nu,R} \leq f_{\nu}$.
Hence, we have a contradiction and we therefore conclude $\nu \geq \nu_0$.

Now let $R > \bar{x}_n$, $h \in (0,R)$ and note that by the comparison principle
\begin{align*}
M(R) - M(R-h) &\geq V(0, \dots, 0, R) - V(0, \dots, 0, R-h)\\
&= \int_{R-h}^{R} f_{\nu,R}\left(t\right) dt \geq \int_{R-h}^{R} f_{\nu_0,R}\left(t\right) dt > 0.
\end{align*}
Hence $M(R)$ is increasing and 
\begin{align}\label{eq:final-contrad}
\frac{M(R) - M(R-h)}{\int_{R-h}^{R} f_{\nu_0}\left(t\right) dt}
\geq
\frac{\int_{R-h}^{R} f_{\nu_0,R}\left(t\right) dt}{\int_{R-h}^{R} f_{\nu_0}\left(t\right) dt}.
\end{align}
Inequality \eqref{eq:final-contrad} holds for any $R > \bar{x}_n$
independent of $h$.
Taking the limit yields
\begin{align*}
\frac{M'(R)}{f_{\nu_0}\left(R\right)} = \liminf_{h \to 0^+} \frac{M(R) - M(R-h)}{\int_{R-h}^{R} f_{\nu_0}\left(t\right) dt}
\geq
\frac{f_{\nu_0,R}\left(R\right)}{f_{\nu_0}\left(R\right)}
\end{align*}
and thus
$$
\liminf_{R \to \infty}\frac{M'(R)}{f_{\nu_0}\left(R\right)} \geq \liminf_{R \to \infty}\frac{f_{\nu_0,R}\left(R\right)}{f_{\nu_0}\left(R\right)}.
$$
This completes the proof of the theorem. $\hfill \Box$\\


\newpage

\section{Applications to well known equations}
\label{sec:applic}

\setcounter{theorem}{0}
\setcounter{equation}{0}

In this section we apply Theorem \ref{th:Phragmen--Lindelof} to some well known PDEs
for which we can solve the ODE in \eqref{eq:ODE} explicitly -- calculate the limit
$$
\liminf_{R\to\infty}\frac{f_{\nu,R}(R)}{f_{\nu}(R)}
$$
and conclude explicit growth estimates.

When stating corollaries for specific classes of PDEs we would sometimes prefer to infer other types of ``weak'' solutions than viscosity solutions whenever such are more suitable or more commonly used for such equations in the literature.
As the equivalence of different kinds of "weak'' solutions often is a nontrivial problem
we will in some cases avoid going into these details,
but this should not make things unclear.
The reason is that we only use comparison between ``weak" subsolutions and classical strict supersolutions -- i.e. Lemma \ref{le:comp-weak-some-sence}.

We begin with the simplest case $\Phi \equiv 0$, including e.g. the famous $p$-Laplace equation, proceed with PDEs having sublinear growth in the gradient according to $\Phi(t,s) = C(t)s^k$ for $k\geq 1$ and end
by the variable exponent $p$-Laplace equation, which satisfies assumption \eqref{eq:ass_drift_super} with $\Phi(s) = C(t) s |\log s|$.

\subsection*{The case $\Phi(s) \equiv 0$}

In this simple case the ODE  \eqref{eq:ODE} reduces to
$$
\frac{df}{dt} = - K \frac{\Lambda(t)}{\lambda(t)}f(t) \quad \text{with} \quad f(0) = \nu
$$
and hence $f_{\nu,R}(t) = \nu e^{-K \int_{0}^{t} \Lambda(s) \lambda^{-1}(s) ds}$ and $f_{\nu}(t) \equiv \nu$.
We obtain the limit
\begin{align}\label{eq:limit_phi=0}
\liminf_{R \to \infty} \frac{M'(R)}{f_{\nu}(R)} \geq
\liminf_{R \to \infty} \frac{f_{\nu,R}(R)}{f_{\nu}(R)} = \lim_{R \to \infty}  e^{-K\int_{0}^{R} \Lambda(s) \lambda^{-1}(s) ds}
\end{align}
which is positive if
$
K \int_{0}^{R} \Lambda(s) \lambda^{-1}(s) ds \precsim 1.
$
Recalling from the definitions in \eqref{eq:ODE} that $K = \frac{n}{\gamma(R)}$ this forces us to chose the function $\gamma(R)$ in the definition of $D(R)$, given by \eqref{eq:domain_2},
\begin{align}\label{eq:choice-of-gamma-phi=0}
\int_{0}^{R}\frac{\Lambda(s)}{\lambda(s)} ds \precsim \gamma(R).
\end{align}
Following the remark just below Theorem \ref{th:Phragmen--Lindelof} we see that \eqref{eq:limit_phi=0} implies
\begin{align}\label{eq:linear-phragmen2}
\liminf_{R\to\infty} \, \frac{M(R)}{R} \, > \, 0,
\end{align}
and we thus retrieve the classical form of the Phragmen--Lindel\"of theorem.
If the PDE is uniformly elliptic, i.e. $\Lambda/\lambda = constant$,
then according to \eqref{eq:choice-of-gamma-phi=0} we can pick $\gamma(R) = R$ and thereby $D(R)$
preserves its geometric proportions for all $R > 0$, which also agrees with the classical Phragmen--Lindel\"of theorem.
If ellipticity blows up at infinity, i.e. $\Lambda(R)/\lambda(R) \to \infty$, then the loss in estimate \eqref{eq:linear-phragmen2} comes only in
the shape of $D(R)$ -- it expands faster in $x'$-direction since we need to take a larger $\gamma(R)$ according to \eqref{eq:choice-of-gamma-phi=0}.

Concerning sharpness of \eqref{eq:linear-phragmen2} we note that with $\phi \equiv 0$ the solution of \eqref{eq:candidate}
yields
$$
\int_{0}^{x_n} \, f_{\nu}(s)\, ds = \nu\, x_n,
$$
which clearly hits the bottom of \eqref{eq:linear-phragmen2}.

Following \eqref{eq:main} and \eqref{eq:ass_drift_super} we see that \eqref{eq:linear-phragmen2} holds, e.g.,
for subsolutions of the quasilinear equations
\begin{align}\label{eq:nysummaexempel}
-\sum_{i,j = 1}^{n} A_{ij}(x) \frac{\partial^2 u}{\partial x_i \partial x_j} + f(x,u,Du) = 0,
\end{align}
corresponding to $F(x,r,p,X) = -\text{Tr}\left(A(x) X\right) + f(x,u,Du)$,
where $A(x) \in \mathbb{S}^n$ satisfies $\lambda(x_n) \text{Tr}\left(Y\right) \leq \text{Tr}\left(A(x) Y\right) \leq \Lambda(x_n)\text{Tr}\left(Y\right)$ for all $Y \geq 0$, and
\begin{align}\label{eq:pucci-ex}
P^-_{\lambda, \Lambda}(D^2u) + f(x,u,Du) = 0,
\end{align}
whenever $f(x,u,D u) \geq 0$ is nondecreasing in $u$.
One such PDE is the following $p$-Laplace equation, $p \in (1, \infty)$, with lower order terms
\begin{align}\label{eq:p-lap-ex}
- \nabla \cdot \left(|D u|^{p-2} D u\right) + f(x,u,Du) =  0. 
\end{align}
%
Indeed, with
\begin{align*}
-F(x,u,Du,D^2u) = \Delta u + (p - 2)\Delta_{\infty}u - \frac{f(x,u,Du)}{|Du|^{p-2}},
\end{align*}
where $\Delta_{\infty}u = \langle D^2 u \frac{Du}{|Du|}, \frac{Du}{|Du|}\rangle$
denotes the infinity Laplace operator,
%
%
we see that $F$ satisfies \eqref{eq:ass_drift_super_old} with $\phi \equiv 0$ and
$$
\min\left\{1,p-1\right\} \text{Tr}(Y) \leq F(x,u,Du, X) - F(x,u,Du, X + Y) \leq \max\left\{1,p-1\right\} \text{Tr}(Y),
$$
whenever $Y\geq 0$. 
Hence $F$ satisfies \eqref{eq:ass_ellipt} with $\lambda = \min\left\{1,p-1\right\}$ and $\Lambda =  \max\left\{1,p-1\right\}$ and Theorem \ref{th:Phragmen--Lindelof} applies.

Recalling that Lemma \ref{le:comp-weak-some-sence} holds for weak solutions (defined in the usual way) to $p$-Laplace type problems or that viscosity solutions and weak solutions are equivalent for some $p$-Laplace type equations
, see e.g.
Juutinen--Lindqvist--Manfredi \cite{JLM01},
Julin--Juutinen \cite{JJ12} and
Medina--Ochoa \cite{MO17},
we retrieve the well known Phragmen--Lindel\"of results in Lindqvist\cite{L85} in the setting of halfspaces.

\subsection*{The case $\Phi(t,s) = C(t)s^k$}

We now consider equations satisfying \eqref{eq:ass_drift_super} with $\Phi(t,s) = C(t)s^k$ for some $k \geq 1$.
Note that such $\Phi$ satisfies \eqref{eq:ass_phi_near} and
hence we expect that Theorem \ref{th:Phragmen--Lindelof} implies that subsolutions must be increasing on the boundary of $D(R)$.
The ODE in \eqref{eq:ODE} yields
$$
\frac{df}{dt} = - \frac{C(t) f^k}{\lambda(t)} - K \frac{\Lambda(t)}{\lambda(t)} f, \quad t \in (0,R),  \quad \text{with} \quad f(0) = \nu.
$$
As $C(t)$ is nonincreasing (by assumptions on $\Phi$) and $\Lambda$ is nondecreasing we can replace the above equation with the separable ODE
$$
\frac{df}{dt} = - A(t) \left( f^k + \widetilde{K} f \right),
$$
where $A(t) := \lambda^{-1}(t) C(t)$ and $\widetilde{K} := \frac{K \Lambda(R)}{C(R)}$.
This is possible since solutions of this ODE will approach the origin faster when $t \in (0,R)$,
and hence it creates a lower bound on the limit in Theorem \ref{th:Phragmen--Lindelof}.
To find the solution for $k>1$ we observe that
$$
\frac{1}{\widetilde K}\int_{\nu}^{f(t)} \left(\frac{1}{y} - \frac{y^{k-2}}{y^{k-1} + \widetilde K} \right) dy = - \int_{0}^{t} A(s) ds
$$
and
$$
\frac{1}{k-1}\left[\log y^{k-1} - \log\left(y^{k-1} + \widetilde{K}\right) \right]_\nu^{f(t)} = - \widetilde{K}\int_{0}^{t} A(s) ds.
$$
Thus
$$
f_{\nu,R}(t) =  \left\{
\begin{array}{ll}
\nu e^{- \left(1 + \widetilde{K}\right)\int_{0}^{t}A(s)ds} & \text{if $k = 1$}, \\
\widetilde{K}^{\frac1{k-1}}\left({e^{(k-1) \widetilde{K} \int_{0}^{t}A(s)ds} \left(\frac{\widetilde{K}}{\nu^{k-1}} + 1\right) - 1} \right)^\frac{1}{1-k} & \text{if $k > 1$},
\end{array}
\right.
$$
and by solving $\frac{df}{dt} = - \lambda(t)^{-1}C(t) f^k$ with $f(0) = \nu$ we also obtain
\begin{align}\label{eq:deriv-sublinear}
f_{\nu}(t) =  \left\{
\begin{array}{ll}
\nu e^{-\int_{0}^{t}A(s)ds } & \text{if $k = 1$}, \\
\left((k - 1) \int_{0}^{t}A(s)ds + \nu^{1-k}\right)^\frac{1}{1-k} & \text{if $k > 1$}.
\end{array}
\right.
\end{align}
The limit in Theorem \ref{th:Phragmen--Lindelof} becomes, for $k=1$,
\begin{align*}
\liminf_{R \to \infty} \frac{M'(R)}{f_{\nu}(R)} \geq
\liminf_{R \to \infty} \frac{f_{\nu,R}(R)}{f_{\nu}(R)} = \lim_{R \to \infty} e^{-\widetilde{K} \int_0^R A(s) ds},
\end{align*}
and for $k > 1$,
\begin{align*}
\liminf_{R \to \infty} \frac{M'(R)}{f_{\nu}(R)} \geq
\liminf_{R \to \infty} \frac{f_{\nu,R}(R)}{f_{\nu}(R)} 
 = \lim_{R \to \infty} \left( \frac{\widetilde{K}(k-1)\int_{0}^{R}A(s)ds + \widetilde{K}\nu^{1-k}}{e^{\widetilde{K}(k-1)\int_{0}^{R}A(s)ds}\left(\widetilde{K}\nu^{1-k} + 1\right) - 1} \right)^{\frac{1}{k-1}}.
\end{align*}
Let's observe that if
%
$
\widetilde{K}(R) \int_{0}^{R}A(s)ds \precsim 1,
$
which is obtained by taking $\frac{\Lambda(R)}{C(R)} \int_{0}^{R}A(s)ds  \precsim \gamma(R)$,
then the limits are positive and we obtain
\begin{align}\label{eq:conclusion-sublinear}
\liminf_{R\to\infty}\frac{M'(R)}{f_\nu(R)} > 0.
\end{align}
Thus, we may derive several Phragmen--Lindel\"of type results using Theorem \ref{th:Phragmen--Lindelof}, whose form will depend on the exponent
$k$ and the functions $C$, 
$\lambda$ and $\Lambda$.
For example, using \eqref{eq:deriv-sublinear}-\eqref{eq:conclusion-sublinear} we have proved:
\begin{corollary}
\label{cor:sublinear1}
Suppose that \eqref{eq:ass_drift_super} holds with $\Phi(t,s) = C(t) s^k$, $k\geq1$.
Let $u$ be a subsolution of \eqref{eq:main} in $\mathbb{R}^n_+$ satisfying
\begin{align*}
\limsup_{x \to y} \, u(x) \, \leq \,  0 \quad \textrm{for all} \quad y \in \partial \mathbb{R}^n_+. 
\end{align*}
Assume also that $u(\bar{x}) > 0$ for some $\bar{x}$ on the $x_n$-axis.
Then the following is true, with $A(t) = C(t)\lambda(t)^{-1}$:
\begin{itemize}
\item[$(i)$] 
If $\int_{0}^{R}A(s)ds \precsim R^{\alpha (k-1)}$ for $\alpha \geq 0, k>1$ and $\frac{\Lambda(R)}{\gamma(R)C(R)}\precsim R^{-\alpha(k-1)}$ then
\begin{align*}
\liminf_{R\to\infty} \, \frac{M'(R)}{R^{-\alpha}} \, > \, 0 \quad \text{implying} \quad
\liminf_{R\to\infty} \, \frac{M(R)}{R^{1-\alpha}} \,>\,0.
\end{align*}
\item[$(ii)$] If $\int_{0}^{R}A(s)ds \precsim R$ and $\frac{\Lambda(R)}{\gamma(R)C(R)}\precsim\frac1R$ then
%
\begin{align*}
\liminf_{R\to\infty} \, \frac{M'(R)}{e^{-R}} \, > \, 0 \quad \text{if $k = 1$}, \quad 
\liminf_{R\to\infty} \, \frac{M'(R)}{R^{\frac{1}{1-k}}} \, > \, 0 \quad \text{if $k\in(1,2)$},
\end{align*}
%
%
\begin{align*}
\liminf_{R\to\infty} \, \frac{M(R)}{\log(R)} \, > \, 0 \quad \text{if $k=2$}
\quad \text{and} \quad
\liminf_{R\to\infty} \, \frac{M(R)}{R^{\frac{2-k}{1-k}}} \, > \, 0 \quad \text{if $2 < k$}.
\end{align*}
\item[$(iii)$] 
 If $\int_{0}^{R}A(s)ds \precsim \log(R)$, $\frac{\Lambda(R)}{\gamma(R)C(R)} \precsim\frac1{\log(R)}$ and $k=1$ then
\begin{align*}
\liminf_{R\to\infty} \, \frac{M(R)}{\log(R)} \, > \, 0.
\end{align*}
\end{itemize}
\end{corollary}

We remark that in Corollary \ref{cor:sublinear1} we only summarize some examples of growth estimates that take simple forms -- the reader may return to conclusion \eqref{eq:conclusion-sublinear} for the general case.
Note also that all conclusions in Corollary \ref{cor:sublinear1} are independent of $\nu$,
meaning that we only have to use arbitrary small $\nu > 0$ to prove them.
Therefore, since $f_{\nu}$ and $f_{\nu,R}$ are nonincreasing functions in this case,
we only need that the assumption $\Phi(t,s) = C(t)s^k$ holds for arbitrary small $s$.
%

Conclusion $(i)$ takes the form of the classical
Phragmen--Lindel\"of theorem and when $\alpha = 0$ it applies e.g. when
\begin{align}\label{eq:tjohej-sublinear}
\Phi(|x|,s) = C(|x|) s^k = \frac{c}{(1+|x|)^a} s^k,
\end{align}
$k \geq 1$, $a > 1$,
ellipticity $\lambda(R) = constant$, $\Lambda(R) = constant$ and $\gamma(R) \succsim R^a$.
Conclusion $(ii)$ holds e.g. when $A = C/\lambda = constant$, $\Lambda/C = constant$ and $\gamma(R) \equiv R$.
We observe that the exponent $k$ in $\Phi(t,s) = C(t) s^k$ has a borderline value at $k=2$.
Namely, if $k \in [1,2)$ then subsolutions may be bounded, but if $k\in [2,\infty)$ then any subsolution must grow to infinity.
As already mentioned in Section \ref{sec:res} such border is,
beyond the assumptions in Corollary \ref{cor:sublinear1},
characterized by convergence/divergence of
\begin{align*}
\int_{0}^{\infty} f_{\nu}(s) ds.
\end{align*}
Conclusion $(iii)$ holds e.g. when $a = 1$ in \eqref{eq:tjohej-sublinear},
$\gamma(R) = R\log(R)$ and $\Lambda = constant$.

We further remark that upper bounds on $\int_{0}^{R}A(s)ds$ have played an important role
for related results in the literature, see e.g. Gilbarg \cite{G52}, Hopf \cite{H52} and Vitolo \cite{V04}, and that Phragmen--Lindel\"of theorems for similar equations in more general domains but with $k=1$ and $k=2$ are proved by Capuzzo-Dolcetta--Vitolo \cite{CDV07} and Koike--Nakagawa \cite{KN09}.
%

As in the case $\Phi \equiv 0$ the results in this subsection apply to PDEs of type \eqref{eq:nysummaexempel} and \eqref{eq:pucci-ex} but now with relaxed assumption on $f$, namely
$
f(x,u,Du) \geq -C(|x|)|Du|^{k}.
$
In case of the $p$-Laplace type equation \eqref{eq:p-lap-ex}, $1 < p < \infty$,
the growth condition on the lower order terms becomes $f(x,u,Du) \geq -C(|x|)|Du|^{\,k+p-2}$.

\begin{figure}[!hbt]
\begin{center}
\includegraphics[height = 7cm, width = 7 cm]{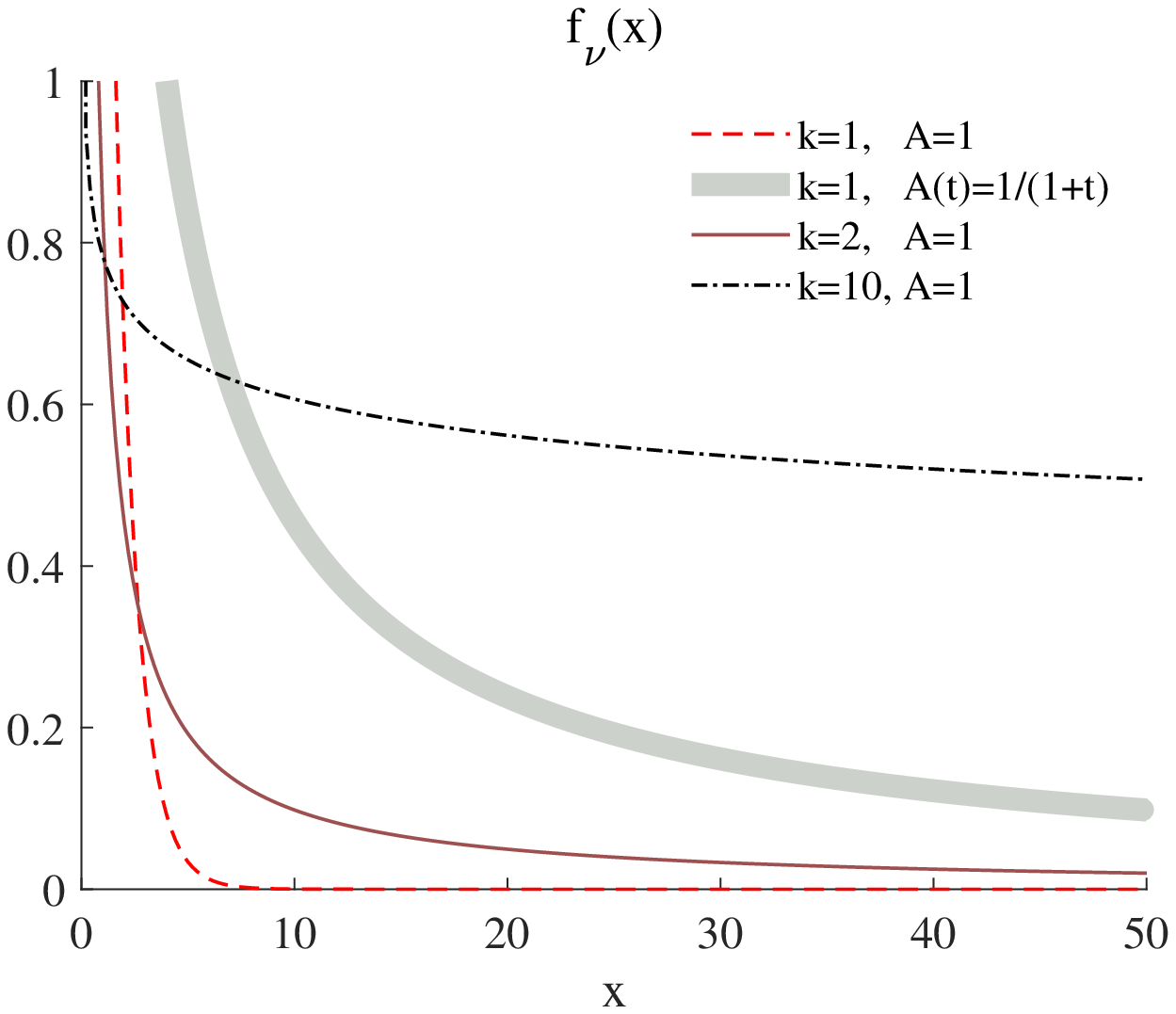}
\hspace{0.8cm}
\includegraphics[height = 7cm, width = 7 cm]{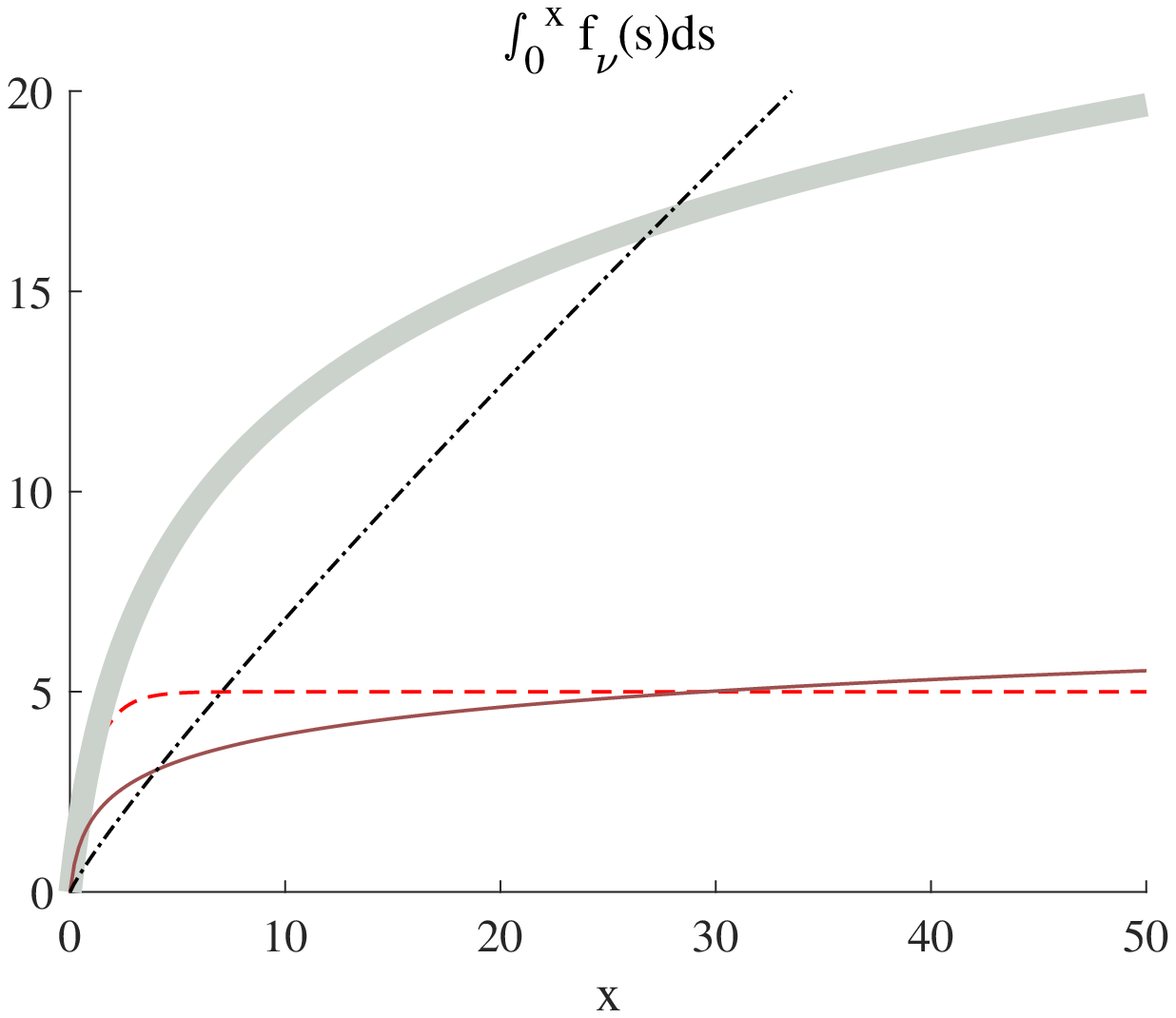}
\end{center}
\caption{The derivative $f_\nu$ in \eqref{eq:deriv-sublinear} and the solution $\int_0^{x} f_{\nu}(t) dt$. Solid curves (for $u$) approach infinity at speed $\log(x)$ and dashed-dot at speed $x^{8/9}$, while dashed curve is bounded.
Corresponding growth estimates are established in Corollary \ref{cor:sublinear1}. In all simulations, $\nu = 5.$}
\label{fig:s^k}
\end{figure}

When $A(s) = C(s)/\lambda(s) = constant$ we can explicitly find the classical solution of \eqref{eq:candidatesPDE}, ensuring sharpness, in case of $\Phi(t,s) = C(t)s^k$:
\begin{align*}
u(x_n) = \int_0^{x_n} f_{\nu}(t) dt = \frac{1}{A} \left\{
\begin{array}{ll}
 \nu\left(1 - e^{-A x_n}\right)  & \text{if $k = 1$},\\
\log\left( A x_n \nu + 1\right) & \text{if $k = 2$}, \\
\frac{\nu^{2-k}}{2-k} \left( 1 - \left((k-1) A x_n \nu^{k-1} + 1\right)^\frac{2-k}{1-k}\right) & \text{otherwise}.
\end{array}
\right.
\end{align*}
Figure \ref{fig:s^k} shows the solutions $\int_0^{x_n} f_{\nu}(t) dt$ for some values of $k$, $\nu$ and different functions $A(t)$.
%


\subsection*{The case $\Phi(t,s) = C(t) s|\log s|$ : variable exponent $p$-Laplace equation}

We set $\Phi(t,s) = C(t) s|\log s|$ and obtain the ODE
\begin{align*}
\frac{df}{dt} = - \frac{C(t) f|\log f|}{\lambda(t)} - K \frac{\Lambda(t)}{\lambda(t)}  f, \quad t \in (0,R),  \quad \text{with} \quad f(0) = \nu.
\end{align*}
By the same argument as in the case $\Phi = C(t) s^k$ we replace the above ODE with
$$
\frac{df}{dt} = - A(t) \left( f|\log f| - \widetilde{K} f \right),
$$
where $A(t) = \lambda^{-1}(t) C(t)$ and $\widetilde{K} = \frac{K \Lambda(R)}{C(R)}$.
This equation separates, when $0 < f \leq 1$, to
$$
\log\left(\widetilde{K}-\log \nu\right) - \log\left(\widetilde{K}-\log f\right) = - \int_0^t A(s) ds.
$$
Thus
$$
f_{\nu,R}(t) = e^{\widetilde{K}\left(1 - e^{\int_0^t A(s) ds}\right)}\nu^{e^{\int_0^t A(s) ds}}
$$
and since $f$ must be nonincreasing this holds for $0 < \nu \leq 1$.
If $f > 1$ the solution takes a similar form.
In particular
\begin{align}\label{eq:derivative-p(x)-help}
f_{\nu,R}(t) = \left\{
\begin{array}{ll}
e^{\widetilde{K}\left(1 - e^{\int_0^t A(s) ds}\right)} \nu^{e^{\int_0^t A(s) ds}}   & \text{if $0 < \nu \leq 1$},\\
e^{\widetilde{K}\left(e^{-\int_0^t A(s) ds} - 1\right)} \nu^{e^{-\int_0^t A(s) ds}}  & \text{if $1 < f(t)$}
\end{array}
\right.
\end{align}
and
\begin{align}\label{eq:derivative-p(x)}
f_{\nu}(t) = \left\{
\begin{array}{ll}
 \nu^{e^{\int_0^t A(s) ds}}   & \text{if $0 < \nu \leq 1$},\\
 \nu^{e^{-\int_0^t A(s) ds}}  & \text{if $1 < \nu$}.
\end{array}
\right.
\end{align}
The limit in Theorem \ref{th:Phragmen--Lindelof} becomes, for $0 < \nu \leq 1$,
\begin{align}\label{eq:limit p(x)-ny}
\liminf_{R \to \infty} \frac{M'(R)}{f_{\nu}(R)} \geq
\liminf_{R \to \infty} \frac{f_{\nu,R}(R)}{f_{\nu}(R)} =
\lim_{R \to \infty} e^{\widetilde{K}\left(1 - e^{\int_0^R A(s) ds}\right)}
\end{align}
which is positive if $\widetilde{K} e^{\int_0^R A(s) ds} \precsim 1$.
Therefore, we have to pick $\gamma(R) \succsim \frac{\Lambda(R)}{C(R)} e^{\int_0^R A(s) ds}$ to achieve a growth estimate.

When $\nu > 1$ we know that $f_{\nu,R}$ in \eqref{eq:derivative-p(x)-help} stays above 1 if
$$
\left( \widetilde{K} + \log \nu\right) e^{-\int_0^t A(s) ds} > \widetilde{K}
$$
which forces us to take $\gamma(R) \succsim \frac{\Lambda(R)}{C(R)} e^{\int_0^R A(s) ds}$.
Then
\begin{align*}
\liminf_{R \to \infty} \frac{M'(R)}{f_{\nu}(R)} \geq
\liminf_{R \to \infty} \frac{f_{\nu,R}(R)}{f_{\nu}(R)} =
\lim_{R \to \infty} e^{\widetilde{K}\left(e^{-\int_0^R A(s) ds} - 1\right)} \geq
\lim_{R \to \infty} e^{-\widetilde{K}} > 0
\end{align*}
as $\gamma(R) \succsim \frac{\Lambda(R)}{C(R)}$.
If the solution $f_{\nu,R}$ decreases to 1 then the limit can be estimated as in \eqref{eq:limit p(x)-ny} since $f_{\nu,R}$ follows the expression for $\nu \in (0,1]$ with $\nu = 1$.
%
%
%
%


In summary, since $\Phi(t,s) = C(t) s |\log s|$ satisfies \eqref{eq:ass_phi_near}
we can conclude that for a subsolution $u$ satisfying the assumptions in Theorem \ref{th:Phragmen--Lindelof} with $\Phi(t,s) = C(t) s |\log s|$,
the following is true when $\gamma(R) \succsim \frac{\Lambda(R)}{C(R)} e^{\int_0^R A(s) ds}$:

\begin{itemize}
\item If $u(\bar x) \geq \check{u}$ for some $\bar{x}$ on the $x_n$-axis and $\nu \in (0,1]$, 
then $M(R)$ may be bounded but
\begin{align}\label{eq:tlogt-1}
\liminf_{R\to\infty} \, \frac{M'(R)}{\nu^{e^{\int_0^R A(s) ds}}} > 0.
\end{align}
\item If $u(\bar x) \geq \check u$ for some $\bar{x}$ on the $x_n$-axis and $\nu  > 1$, 
then $M(R)$ approaches infinity according to
\begin{align}\label{eq:tlogt-2}
\liminf_{R\to\infty} \, \frac{M'(R)}{\nu^{e^{-\int_0^R A(s) ds}}} > 0 \quad \text{implying} \quad
\liminf_{R\to\infty} \, \frac{M(R)}{R} > 0.
\end{align}
\end{itemize}
Thus we retrieve the classical form of a Phragmen--Lindel\"of theorem if the subsolution exceeds
$\check{u}$ in \eqref{eq:solution-variable-p} with $\nu \geq 1$,
but if the subsolution only exceeds $\check{u}$ with  $\nu < 1$, it may grow very slowly
and need not approach infinity.
The border at $\nu = 1$ originates from the fact that $\Phi(1) = 0$ and thus $f_{1} \equiv 1$,
while $\Phi(s) > 0$ for all other positive $s$ implying $f_{\nu} \to 0$ if $\nu \in (0,1)$ and $f_{\nu} \to 1$ if $\nu > 1$.

The solution of \eqref{eq:candidatesPDE} with $\Phi(t,s) = C(t) s |\log s|$, i.e.
\begin{align}\label{eq:denlillarackarn}
C(t) |Du| | \log |Du| | + \lambda(x_n) \Delta u = 0,
\end{align}
ensuring sharpness, can be calculated analytically when $A(t) = C(t)/\lambda(t)= constant$ and then yields
\begin{align}\label{eq:solution-variable-p}
\check u(x) = \int_{0}^{x_n} f_\nu(s) ds = \frac{1}{A} \left\{
\begin{array}{ll}
-E_i\left(\log\nu\right) + E_i\left(e^{A t}\log\nu\right)   & \text{if $0 < \nu < 1$},\\
x_n  & \text{if $\nu = 1$},\\
E_i\left(\log\nu\right) - E_i\left(e^{-At}\log\nu\right)  & \text{if $1 < \nu$}
\end{array}
\right.
\end{align}
where $E_i$ is the Exponential integral.
See Figure \ref{fig:variablep} (upper row) for some illustrations of the functions
$f_\nu$ in \eqref{eq:derivative-p(x)} and $\check{u}$ in \eqref{eq:solution-variable-p}.

\subsubsection*{Variable exponent $p$-Laplace equation}

The $p(x)$-Laplace equation, which
often serves as a model example for PDEs with nonstandard growth, yields
\begin{equation}\label{eq:intro-p(x)}
\nabla\cdot(|D u|^{p(x)-2}D u)\,=\,0.
\end{equation}
The function $p:\Omega\to (1,\infty)$ is usually called a variable exponent.
If $p=constant$, then this equation is the classical $p$-Laplace equation and if $p = 2$ it's the
famous Laplace equation.
Apart from interesting theoretical considerations such equations arise in the applied sciences, for instance
in fluid dynamics, see e.g. Diening--R\r u\v zi\v cka~\cite{dr}, in the study of image processing, see e.g. Chen-Levine-Rao~\cite{clr} and for electro-rheological fluids,
we refer the reader to Harjulehto--H\"ast\"o--L\^e--Nuortio~\cite{hhn} for a recent survey and further references.

We recall the following standard definition of weak solutions of \eqref{eq:intro-p(x)}:
 A function $u \in W^{1, p(x)}_{loc}(\Omega)$ is a weak (sub)solution of \eqref{eq:intro-p(x)} if
\begin{equation*}
\int_\Omega |D u|^{p(x)-2} D u \cdot D \psi\, dx (\leq) = 0
\end{equation*}
for all (nonnegative) $\psi \in C_0^\infty(\Omega)$.
Similarly, $u$ is a weak supersolution if $-u$ is a weak subsolution.
A function which is both a weak subsolution and a weak supersolution is called a weak solution.
An (USC/LSC) weak (sub/super)solution is called a $p(x)$-(sub/super)harmonic function.
We also note that $u \in W^{1,1}_{loc}(\Omega)$ is $p(x)$-harmonic in $\Omega$ if it is a local minimizer of the energy
$$
\int_{\Omega} \frac{1}{p(x)} |Du|^{p(x)} dx,
$$
where $1 < p(x) < \infty$.

To proceed we define the operator
\begin{align*}
\Delta_{p(x)} u := \Delta u + (p(x) - 2)\Delta_{\infty}u + \log |Du| \langle Dp(x), Du\rangle,
\end{align*}
where $\Delta_{\infty}u = \langle D^2 u \frac{Du}{|Du|}, \frac{Du}{|Du|}\rangle$
denotes the infinity Laplace operator.
We note that
%
$
\Delta_{p(x)} u 
\geq 0
$
%
implies
\begin{align}\label{eq:new-explanation}
 -\widehat{F}(x,Du,D^2u) := \Lambda(x) \text{Tr}(D^2u^+) - \lambda(x) \text{Tr}(D^2u^-) + |Dp(x)| |Du| |\log |Du|| \geq 0
\end{align}
with
$\lambda(x) = \min\{1, p(x) - 1\}$ and
$\Lambda(x) = \max\{1, p(x) - 1\}$. 
This suggests that $p(x)$-subharmonic functions should be viscosity subsolutions to $\widehat{F} = 0$, which is the case.
Indeed, following the proof in Julin \cite{J13}, which expands on
Juutinen--Lukkari--Parviainen \cite{JLP10}, we can conclude the following slightly generalized version of \cite[Lemma 5.2]{J13}:


\begin{lemma}\label{te:Julin-generalized}
Suppose that $p(x)$ is %
$C^1(\mathbb{R}_+)$, 
$1 <  p(x) < \infty$,
$\lambda(x) = \min\{1, p(x) - 1\}$ and
$\Lambda(x) = \max\{1, p(x) - 1\}$.
If $u$ is $p(x)$-subharmonic in a domain $\Omega \in \mathbb{R}^n$,
$\varphi \in C^2(\Omega)$ is such that
$\varphi(x_0) = u(x_0)$ at $x_0 \in \Omega$ and $\varphi \geq u$ then
$$
\widehat{F}(x_0,D\varphi(x_0),D^2(\varphi(x_0)))
\leq -\Delta_{p(x)}\varphi(x_0) \leq 0.
$$
%
%
%
\end{lemma}
To obtain a PDE satisfying the required assumptions we redefine $\widehat{F}$ by replacing ellipticity with $\lambda(x_n) \leq \min\{1, p(x) - 1\}$ nonincreasing,
$\Lambda(x_n) \geq \max\{1, p(x) - 1\}$ nondecreasing
and also replacing the nonhomogeneous term with $\Phi(|x|,s) \geq |Dp(x)| s |\log s|$, where $\Phi(|x|,s)$ is nonincreasing in $|x|$.
In particular, we can take
\begin{align}\label{eq:functions-p(x)}
\lambda(t) =& p_\lambda := \inf_{y:y_n\leq t}\min\{1, p(y) - 1\}, \quad \Lambda(t) = p_\Lambda := \sup_{y:y_n\leq t}\max\{1, p(y) - 1\} \quad \text{and} \notag\\
\Phi(t,s) &= ||Dp||_{\infty,t}\, s |\log s|, \quad \text{where} \quad ||f||_{\infty,t} = \sup_{y : |y| \geq t} |f(y)|.
\end{align}
%
%
In conclusion, a weak USC subsolution (a $p(x)$-subharmonic function) to the variable exponent $p$-Laplace equation is a viscosity subsolution of a PDE of type \eqref{eq:main} satisfying  \eqref{eq:ass_drift_super} and \eqref{eq:ass_phi_near}.
We can therefore conclude that deductions \eqref{eq:tlogt-1} and \eqref{eq:tlogt-2} hold for $p(x)$-subharmonic functions whenever $p(x)$ is $C^1(\mathbb{R}_+)$ and
$1 < p(x) < \infty$.

We summarize our findings in the following theorem yielding
Phragmen--Lindel\"of-type results,
of which some are sharp,
for weak solutions to the variable exponent $p$-Laplace equation:

\begin{theorem}
\label{thm:p(x)}
Suppose that $p(x)$ is %
$C^1(\mathbb{R}^n_+)$, 
$1 <  p(x) < \infty$, and
let $u$ be $p(x)$-subharmonic in $\mathbb{R}^n_+$
satisfying
\begin{align*}
\limsup_{x \to y} \, u(x) \, \leq \,  0 \quad \textrm{for all} \quad y \in \partial \mathbb{R}^n_+.
\end{align*}
Then $u$ is a viscosity subsolution of an equation of type \eqref{eq:main} satisfying \eqref{eq:ass_drift_super} and \eqref{eq:ass_phi_near} with $\Phi, \lambda = p_\lambda$ and $\Lambda = p_\Lambda$ as in \eqref{eq:functions-p(x)}.
%
%
Moreover, if $\frac{p_\Lambda(R)}{||Dp||_{\infty,R}} {\exp\left({\int_{0}^{R} \frac{||Dp||_{\infty,s}}{p_\lambda(s)}ds}\right)} \precsim \gamma(R)$ and
$\check u(x) = \int_{0}^{x_n} f_\nu(s) ds$ with $f_{\nu}$ from \eqref{eq:derivative-p(x)} then the following is true:
\begin{itemize}
\item If $u \geq \check{u}$ somewhere on the $x_n$-axis for $\nu \in (0,1]$ 
then
\begin{align*}
\liminf_{R\to\infty} \, \frac{M'(R)}{\nu^{\exp\left({\int_{0}^{R} \frac{||Dp||_{\infty,s}}{p_\lambda(s)}ds}\right)}} > 0.
\end{align*}
%
%
\item If $u \geq \check u$ somewhere on the $x_n$-axis for $\nu  > 1$ 
then
\begin{align*}
\liminf_{R\to\infty} \, \frac{M(R)}{R} > 0.
\end{align*}
\end{itemize}
\end{theorem}
%
%
We thus retrieve the classical form of a Phragmen--Lindel\"of theorem if the subsolution exceeds
$\check{u}$ with $\nu \geq 1$,
in particular if it exceeds $x_n$.
On the other hand, if the subsolution only exceeds $\check{u}$ with $\nu < 1$,
then Theorem \ref{thm:p(x)} states that it may grow very slowly and be bounded.
The sharpness in the case $\nu \geq 1$ follows by observing, e.g., that
\begin{align}\label{eq:sharp-new}
u(x) = c\, x_n \quad \text{is $p(x)$-harmonic with} \quad p(x) = M_0 + \sum_{i = 1}^{n-1} M_i x_i^2,
\end{align}
whenever $c \geq 1, M_0 > 1$ and $M_i$, for $i \in [1,n-1]$, are constants.
It is worth observing that the conclusion
\begin{align*}
\liminf_{R\to\infty} \, \frac{M(R)}{R} > 0
\end{align*}
follows also in the case $\nu \in (0,1]$ if ${\int_{0}^{R} \frac{||Dp||_{\infty,s}}{p_\lambda(s)}ds} \precsim 1$ since then
$\liminf_{R\to\infty} \, M'(R) > 0$.
This holds e.g. if the exponent satisfies $p^- < p(x)$
and $||Dp||_{\infty,s} \precsim s^{-k}$ for some constants $p^-, k > 1$;
a natural conclusion since these assumptions force the equation toward the constant exponent $p$-Laplace equation far away from the origin.

Versions of Theorem \ref{thm:p(x)} are possible to derive from \eqref{eq:tlogt-1} and \eqref{eq:tlogt-2};
e.g., it may be useful to redefine the norm in \eqref{eq:functions-p(x)} as $||f||_{\infty,t} = \sup_{y : y_n = t} |f(y)|$.
Then, if $||Dp(x)||_{\infty,x_n} \neq 0$ for all $x_n > 0$
we may divide \eqref{eq:new-explanation} by $||Dp(x)||_{\infty,x_n}$
and conclude, for $\lambda(x_n) \leq \frac{\min\{1,p(x)-1\}}{||Dp(x)||_{\infty,x_n}}$ nonincreasing and
$\Lambda(x_n) \geq \frac{\max\{1,p(x)-1\}}{||Dp(x)||_{\infty,x_n}}$ nondecreasing,
that Theorem \ref{thm:p(x)} holds with
$\frac{p_\Lambda(R)}{||Dp||_{\infty,R}}$ replaced by $\Lambda(R)$
and $\frac{||Dp||_{\infty,s}}{p_\lambda(s)}$ replaced by $\lambda^{-1}(s)$.
In particular, in the case $\nu \in(0,1]$ the conclusion then reads
\begin{align}\label{eq:ny-lambda-johohej}
\liminf_{R\to\infty} \, \frac{M'(R)}{\nu^{e^{\int_{0}^{R} \lambda^{-1}(s)ds}}} > 0.
\end{align}
We build sharpness 
of this result in Remark \ref{cor:tjohej} below
in which we find a family of exponents for which 
the solution in \eqref{eq:solution-variable-p},
which satisfies $M'(R) = \nu^{e^{\int_{0}^{R} \lambda^{-1}(s)ds}}$,
is $p$-harmonic.

%
We further remark that Theorem \ref{thm:p(x)} sharpens some results of Adamowicz \cite{A14} in the geometric setting of halfspaces,
and the $C^1$-assumption on $p(x)$ should be replaceable with locally Lipschitz continuity by approximation arguments.
Furthermore, the reader may recall the remarks made below deductions \eqref{eq:tlogt-1} and \eqref{eq:tlogt-2} and also note that contrary to the results in the former subsection,
for $\Phi(t,s) = C(t) s^k$,
the growth estimates here depend on $\nu$.

Our estimates may not be optimal when
$\log |Du| \langle Dp(x), Du\rangle$ is negative since then we lose information by our choice of $\Phi$.
We can improve by taking $\phi(s) \equiv 0$, but we still lose information when
subsolutions gradients are not ``close to perpendicular" to $Dp(x)$. 
This motivates us to derive better estimates under assumptions excluding e.g. the solution in \eqref{eq:sharp-new}.
We do so by studying a nonpositive $\Phi$; the case $\Phi(s) = - s |\log s|$, in the next section.

We proceed by proving the following result, in which we find the ``slowest growing" $p(x)$-harmonic function, for a given ellipticity bound, and the corresponding family of exponents. 

\begin{remark}
\label{cor:tjohej}
The function $\check u(x) = \int_{0}^{x_n} f_\nu(s) ds$ with $f_{\nu}$ from \eqref{eq:derivative-p(x)} is $p(x)$-harmonic with exponent
$$
\check p(x) = 1 + Me^{-\int_{0}^{x_n}A(s)ds} \quad \text{if $\nu \in (0,1]$ and} \quad
\check p(x) = 1 + Me^{\int_{0}^{x_n}A(s)ds}\quad \text{if $\nu >1$},
$$
%
whenever $M \in \mathbb{R}$ is a constant.

Suppose that $\nu \in (0,1]$ and $A^{-1}(s) = \lambda(s)$ is nonincreasing. Then 
$\check u(x)$ is the slowest growing $p(x)$-harmonic function
in the sense of version \eqref{eq:ny-lambda-johohej} of Theorem \ref{thm:p(x)}.
In particular, any $p(x)$-subharmonic function with exponent
$p(x) \in C^1(\mathbb{R}^n_+)$,
$1 < p(x) < \infty$,
$\lambda(x_n) \leq \frac{\min\{1, p(x) - 1\}}{||Dp||_{\infty, x}}$ and $\frac{\max\{1, p(x) - 1\}}{||Dp||_{\infty, x}} \leq \Lambda(x_n)$, satisfying
\begin{align*}
\limsup_{x \to y} \, u(x) \, \leq \,  0 \quad \textrm{for all} \quad y \in \partial \mathbb{R}^n_+ 
\end{align*}
that exceeds $\check u(x)$ somewhere on the $x_n$-axis satisfies
\begin{align*}
\liminf_{R\to\infty} \, \frac{M'(R)}{\check{u}(R e_n)} = \liminf_{R\to\infty} \, \frac{M'(R)}{\nu^{e^{\int_{0}^{R} \lambda^{-1}(s)ds}}} > 0.
\end{align*}


Finally, if $\lambda$ is constant then
\begin{align*}
\check u(x) = \lambda \left\{
\begin{array}{ll}
-E_i\left(\log\nu\right) + E_i\left(e^{\lambda^{-1}x_n}\log\nu\right)   & \text{if $0 < \nu < 1$},\\
x_n  & \text{if $\nu = 1$}
\end{array}
\right.
\end{align*}
where $E_i$ is the Exponential integral.
\end{remark}

\noindent
{\bf Proof.}
%
Since $\check{u}$ depends only on $x_n$ and solves \eqref{eq:denlillarackarn} the first statement
follows if we prove that the variable exponent $p(x)$-Laplace equation,
with exponent $\check p(x) = 1 + Me^{\mp\int_{0}^{x_n}A(s)ds}$,
reduces to the PDE \eqref{eq:denlillarackarn} in one dimension.
Without derivatives in $x'$-direction we have
\begin{align}\label{eq:variable-p(x)-one-dimension}
\Delta_{p(x)} u(x) = (p(x) - 1)u''_{x_n x_n}(x) + \log |u_{x_n}'(x)|  p'_{x_n}(x) u_{x_n}'(x) = 0.
\end{align}
Observe that exponent $\check{p}(x)$ is the unique family of $C^1(\mathbb{R}_+^n)$ solutions to the ODE
$$
p_{x_n}'(x) = \mp (p(x) - 1)A(x_n)
$$
and substituting this equality into \eqref{eq:variable-p(x)-one-dimension} yields
%
$
 u''_{x_n x_n}(x)  \mp  A(x_n) \log |u_{x_n}'(x)|  u_{x_n}'(x) = 0
$
%
where the ``$-$" sign is for $\nu \in (0,1]$ when $\log |\check u_{x_n}'(x)| < 0$. Thus
\begin{align*}
  u''_{x_n x_n}(x) + A(x_n) |\log |u_{x_n}'(x)||  u_{x_n}'(x) = 0
\end{align*}
which is \eqref{eq:denlillarackarn} in one dimension.

To prove the second statement we need to ensure that a weak subsolution 
of the variable exponent $\check p$-Laplace equation,
$\check p(x) = 1 + Me^{- \int_{0}^{x_n}\lambda(s)^{-1}ds}$ for some $M$,
of which $\check u$ is a solution,
is a viscosity subsolution of \eqref{eq:main} where \eqref{eq:ass_drift_super} holds with
the same $\Phi(s)$ and $\lambda(t)$ as in version \eqref{eq:ny-lambda-johohej} of Theorem \ref{thm:p(x)}.
To do so we observe that, recalling Lemma \ref{te:Julin-generalized},
any $p(x)$-subharmonic function is viscosity solution of
\begin{align*}
\Delta_{p(x)}u =  \Delta u + (p(x) - 2)\Delta_{\infty}u + \log |Du| \langle Dp(x), Du\rangle \geq 0,
\end{align*}
and hence of
\begin{align*}
|Dp||Du| |\log |Du|| + \max\{1,p(x)-1\} Tr((D^2u)^+) - \min\{1,p(x)-1\} Tr((D^2u)^-) \geq 0.
\end{align*}
Inserting $\check p(x) = 1 + Me^{- \int_{0}^{x_n}\lambda(s)^{-1}ds}$, $|D\check p| = (\check p(x)-1)\lambda^{-1}(x_n)$ and assuming that
$1 < \check p(x) \leq 2$, which we may by taking $M \in (0, 1]$, we see that
\begin{align*}
 |Du| |\log |Du|| + \frac{\lambda(x_n)}{M e^{-\int_{0}^{x_n}\lambda^{-1}(s)ds}} Tr((D^2u)^+) -  \lambda(x_n) Tr((D^2u)^-) \geq 0.
\end{align*}
This is a PDE satisfying \eqref{eq:ass_drift_super} with $\Phi(s) = s|\log s|$ and $\lambda(t)$ as in version \eqref{eq:ny-lambda-johohej} of Theorem \ref{thm:p(x)}.

It remains to show that $\check{p}$ satisfies
$$
\lambda(x_n) \leq \frac{\min\{1, \check p(x) - 1\}}{||D\check p||_{\infty, x}}.
$$
This holds with equality since
$$
\check p(x) - 1 = Me^{-\int_{0}^{x_n}\lambda^{-1}(s)ds}, \qquad 
||D\check p||_{\infty, x} = |\check{p}'(x)| = - \lambda(x_n)^{-1}M e^{-\int_{0}^{x_n}\lambda^{-1}(s)ds}
$$
and we have assumed $M \in (0,1]$.
The proof is complete. $\hfill \Box$\\


\subsection*{The case $\Phi(s) = - s|\log s|$ }

In this case the ODE \eqref{eq:ODE} becomes (we skip the $t$-dependence in $\Phi$ only for simplicity)
$$
\frac{df}{dt} =  \frac{f|\log f|}{\Lambda(t)}  - K \frac{\Lambda(t)}{\lambda(t)}  f, \quad t \in (0,R),  \quad \text{with} \quad f(0) = \nu.
$$
By the same argument as in the former cases we replace this ODE by
$$
\frac{df}{dt} =  \Lambda^{-1}(t)\left( f|\log f|  - \widehat K  f \right),
$$
which separates, and we obtain, with $\widehat K = K \frac{ \Lambda^2(R)}{\lambda(R)}$
\begin{align}\label{eq:derivative-p(x)-help-neg}
f_{\nu,R}(t) = \left\{
\begin{array}{ll}
e^{-\widehat K(1 - e^{-\int_{0}^{t}\Lambda^{-1}(s) ds})} \nu^{e^{-\int_{0}^{t}\Lambda^{-1}(s) ds}} & \text{if $0 < \nu \leq 1$},\\
e^{-\widehat K\left(e^{\int_{0}^{t}\Lambda^{-1}(s) ds} - 1\right)} \nu^{e^{\int_{0}^{t}\Lambda^{-1}(s) ds}}  & \text{if $1 < f(t)$}
\end{array}
\right.
\end{align}
and
\begin{align}\label{eq:derivative-p(x)-neg}
f_{\nu}(t) = \left\{
\begin{array}{ll}
 \nu^{e^{-\int_{0}^{t}\Lambda^{-1}(s) ds}}   & \text{if $0 < \nu \leq 1$},\\
 \nu^{e^{\int_{0}^{t}\Lambda^{-1}(s) ds}}  & \text{if $1 < \nu$}.
\end{array}
\right.
\end{align}
The limits become, for $0 < \nu \leq 1$,
\begin{align*}
\liminf_{R \to \infty} \frac{M'(R)}{f_{\nu}(R)} \geq
\liminf_{R \to \infty} \frac{f_{\nu,R}(R)}{f_{\nu}(R)} \geq
\lim_{R \to \infty} e^{-\widehat K}
\end{align*}
and we only need $\widehat K \precsim 1$. 
When $\nu > 1$ we know that $f_{\nu,R}$ in \eqref{eq:derivative-p(x)-help-neg} stays above 1 if
$$
\left(\log \nu - \widehat K \right) e^{\int_{0}^{t}\Lambda^{-1}(s) ds} > - \widehat K
$$
which forces us to take $\widehat K < \log \nu$.
Then
\begin{align*}
\liminf_{R \to \infty} \frac{M'(R)}{f_{\nu}(R)} \geq
\liminf_{R \to \infty} \frac{f_{\nu,R}(R)}{f_{\nu}(R)} =
\lim_{R \to \infty} e^{-\widehat K\left(e^{\int_{0}^{R}\Lambda^{-1}(s) ds} - 1\right)}
\end{align*}
and we need also $\widehat K e^{\int_{0}^{R}\Lambda^{-1}(s) ds} \precsim 1$ to achieve a growth estimate.

We have defined $\widehat K(R) = \frac{n\Lambda(R)^2}{\lambda(R)\gamma(R)}$ in this case.
However, from \eqref{eq:derivative-p(x)-help-neg} it can be shown that $f_{\nu,R}(t)$ is nondecreasing if $\widehat K \leq |\log(\nu)|$.
This means that
$$
\frac{df}{dt} = \Lambda^{-1}(t)\left( f|\log f|  - \widehat K  f \right) \geq 0, \qquad t \in (0,R).
$$
Now, we let $f_{\nu,R}$ solve this ODE in place of \eqref{eq:ODE} and in the proof of Lemma \ref{le:barrier_super} we replace \eqref{eq:tracetjoho} with
\begin{align*}
\text{Tr}\left( \left(D^2 V \right)^+ \right) &= -\Lambda(\Xi)^{-1} \Phi\left(\Xi,f\left(\Xi\right)\right) - \Lambda(\Xi)^{-1}\widehat K  f(\Xi) + \frac{n-1}{|(x', x_n+\gamma)|} f(\Xi),\\ 
\text{Tr}\left( \left(D^2 V \right)^- \right) &= 0.
\end{align*}
By tracing the remaining part of the proof of Lemma \ref{le:barrier_super} we realize that it is enough to pick
$\widehat K(R) = \frac{n\Lambda(R)}{\gamma(R)}$.

As in the former situation the solution of \eqref{eq:candidatesPDE} with $\Phi(s) = -  s |\log s|$
can be calculated analytically when $\Lambda(t) = constant$:
\begin{align}\label{eq:solution-variable-p-neg}
\check u(x) = \int_{0}^{x_n} f_\nu(s) ds = \Lambda \left\{
\begin{array}{ll}
E_i\left(\log\nu\right) - E_i\left(e^{-\Lambda^{-1}t}\log\nu\right)   & \text{if $0 < \nu < 1$},\\
x_n  & \text{if $\nu = 1$},\\
-E_i\left(\log\nu\right) + E_i\left(e^{\Lambda^{-1}t}\log\nu\right)  & \text{if $1 < \nu$}.
\end{array}
\right.
\end{align}
See Figure \ref{fig:variablep} (lower row) for functions
$f_\nu$ in \eqref{eq:derivative-p(x)-neg} and $\check{u}$ in \eqref{eq:solution-variable-p-neg}.

Now, using the calculations above \eqref{eq:new-explanation} we see that $\Delta_{p(x) u \geq 0}$ implies
\begin{align*}
 \max\{1,p(x)-1\} \text{Tr}(D^2u^+) - \min\{1,p(x)-1\} \text{Tr}(D^2u^-) + |Dp|  |Du|  \cos\theta \log |Du| \geq 0
\end{align*}
where $\theta = \theta(x)$ is the angle between $Du$ and $Dp$.
Assume $|Dp| |\cos \theta| > 0$ and divide the PDE with this factor to obtain
\begin{align*}
 \Lambda(x_n) \text{Tr}(D^2u^+) - \lambda(x_n) \text{Tr}(D^2u^-) + \frac{\cos\theta}{|\cos\theta|} |Du|  \log |Du| \geq 0
\end{align*}
where $\lambda(x_n) \leq \frac{\min\{1,p(x)-1\}}{|Dp| |\cos \theta|}$ and $\Lambda(x_n) \geq \frac{\max\{1,p(x)-1\}}{|Dp| |\cos \theta|}$ for some nonincreasing function $\lambda$ and nondecreasing function $\Lambda$.
Assuming $\cos\theta \log|Du| \leq 0$ leads to
\begin{align*}
 \Lambda(x_n) \text{Tr}(D^2u^+) - \lambda(x_n) \text{Tr}(D^2u^-) - |Du|  |\log |Du|| \geq 0
\end{align*}
and we can apply the results from
\eqref{eq:derivative-p(x)-neg}-
\eqref{eq:solution-variable-p-neg} and Lemma \ref{te:Julin-generalized} to obtain:

\begin{corollary}\label{cor:hejehjsista}
Suppose that $p(x)$ and $u$ are as in Theorem \ref{thm:p(x)}.
Let $\theta = \theta(x)$ be the angel between $Dp$ and $Du$ and assume that
$$
|Dp| |\cos \theta| > 0 \qquad \text{and} \qquad \cos\theta \log|Du| \leq 0
$$
holds in $\mathbb{R}^n_+$ (in a suitable weak sense if $u$ is not $C^1$ with $|Du|\neq 0$). 
Then $u$ is a subsolution of an equation of type \eqref{eq:main} satisfying \eqref{eq:ass_drift_super} with
$\Phi(s) =  - s |\log s|$, 
$\lambda(x_n) \leq \frac{\min\{1, p(x) - 1\}}{|Dp| |\cos \theta|}$ and
$\Lambda(x_n) \geq \frac{\max\{1, p(x) - 1\}}{|Dp| |\cos \theta|}$
for some nonincreasing function $\lambda$ and nondecreasing function $\Lambda$.
If $\frac{n \Lambda(R)}{|\log \nu|} < \gamma(R)$ and $\check u(x) = \int_{0}^{x_n} f_\nu(s) ds$ with $f_{\nu}$ from \eqref{eq:derivative-p(x)-neg} then the following is true:
\begin{itemize}
\item If $u \geq \check{u}$ somewhere on the $x_n$-axis, $\nu \in (0,1)$ 
then
\begin{align*}
\liminf_{R\to\infty} \, \frac{M(R)}{R} > 0.
\end{align*}
\item If $u \geq \check{u}$ somewhere on the $x_n$-axis, $\nu  > 1$ and $\Lambda(R) e^{\int_{0}^{R}\Lambda^{-1}(s) ds} \precsim \gamma(R)$, then
\begin{align*}
\liminf_{R\to\infty} \, \frac{M'(R)}{\nu^{e^{\int_{0}^{R}\Lambda^{-1}(s)ds}}} > 0
\end{align*}
which, if $\Lambda = constant$ and $E_i$ is the Exponential integral, implies
\begin{align*}
\liminf_{R\to\infty} \, \frac{M(R)}{E_i\left(e^{\Lambda^{-1} R}\log\nu\right)-E_i\left(\log\nu\right)} > 0.
\end{align*}
\end{itemize}
\end{corollary}

In the one dimensional case Corollary \ref{cor:hejehjsista} shows that if we know that the exponent $p(x)$ is increasing, $p' > 0$,
and that the subsolution satisfies $0 < u' < 1$, then $\liminf_{R\to \infty} u(R)/R > 0$,
which is much stronger than the growth estimates that can be derived from Theorem \ref{thm:p(x)}.
Similarly, if we know that $p' < 0$ and $1 < u'$ then $\liminf_{R\to \infty} u(R)/ {\nu^{e^{\int_{0}^{R}\Lambda^{-1}(s)ds}}}> 0$.
These improvements can be visualized by comparing the right panels in Figure \label{fig:variablep};
the upper right panel corresponds to Theorem \ref{thm:p(x)} while the lower right panel corresponds to the results in Corollary \ref{cor:hejehjsista}.

\begin{figure}[!hbt]
\begin{center}
\includegraphics[height = 6.5cm, width = 7 cm]{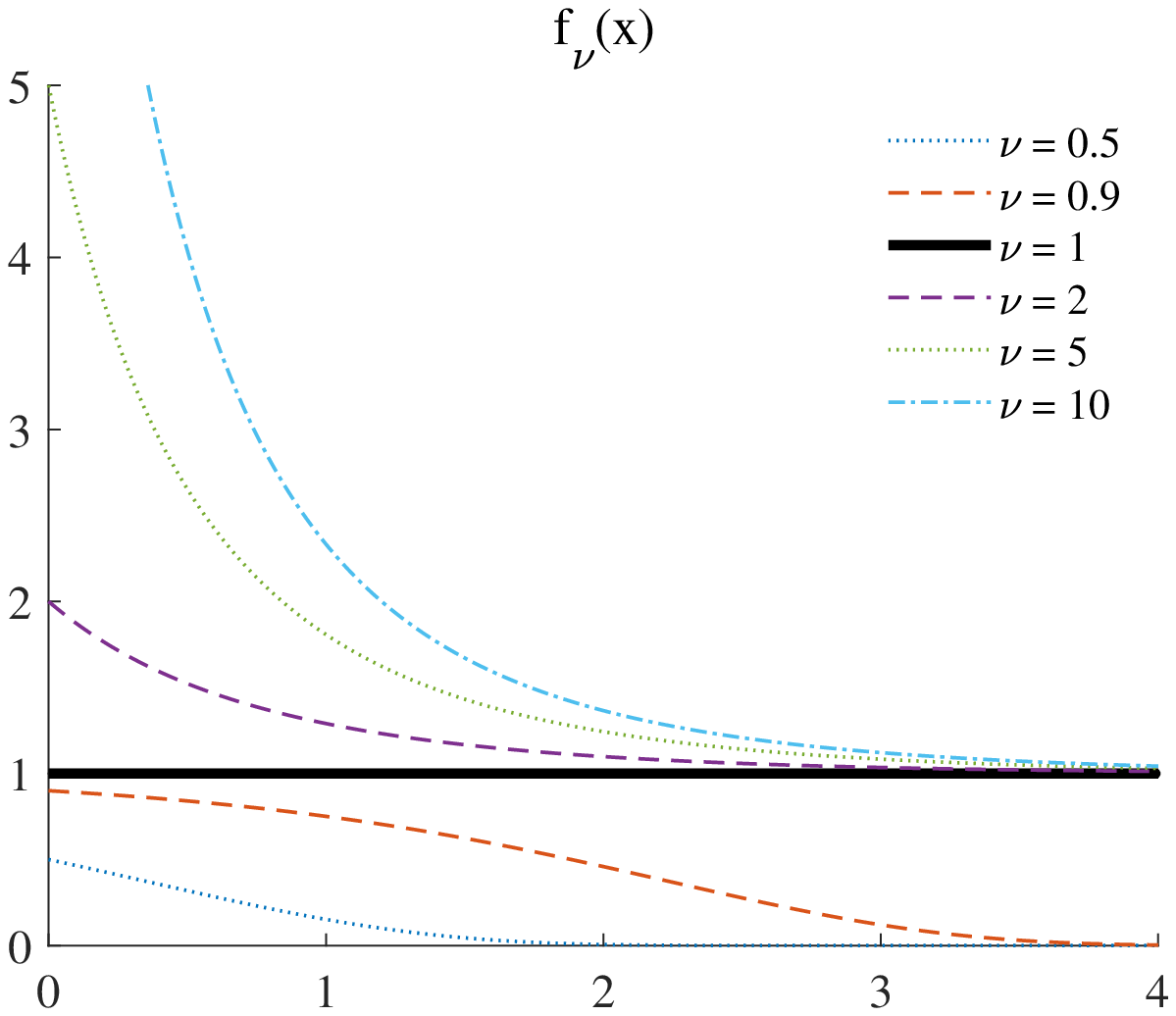}
\hspace{0.8cm}
\includegraphics[height = 6.5cm, width = 7 cm]{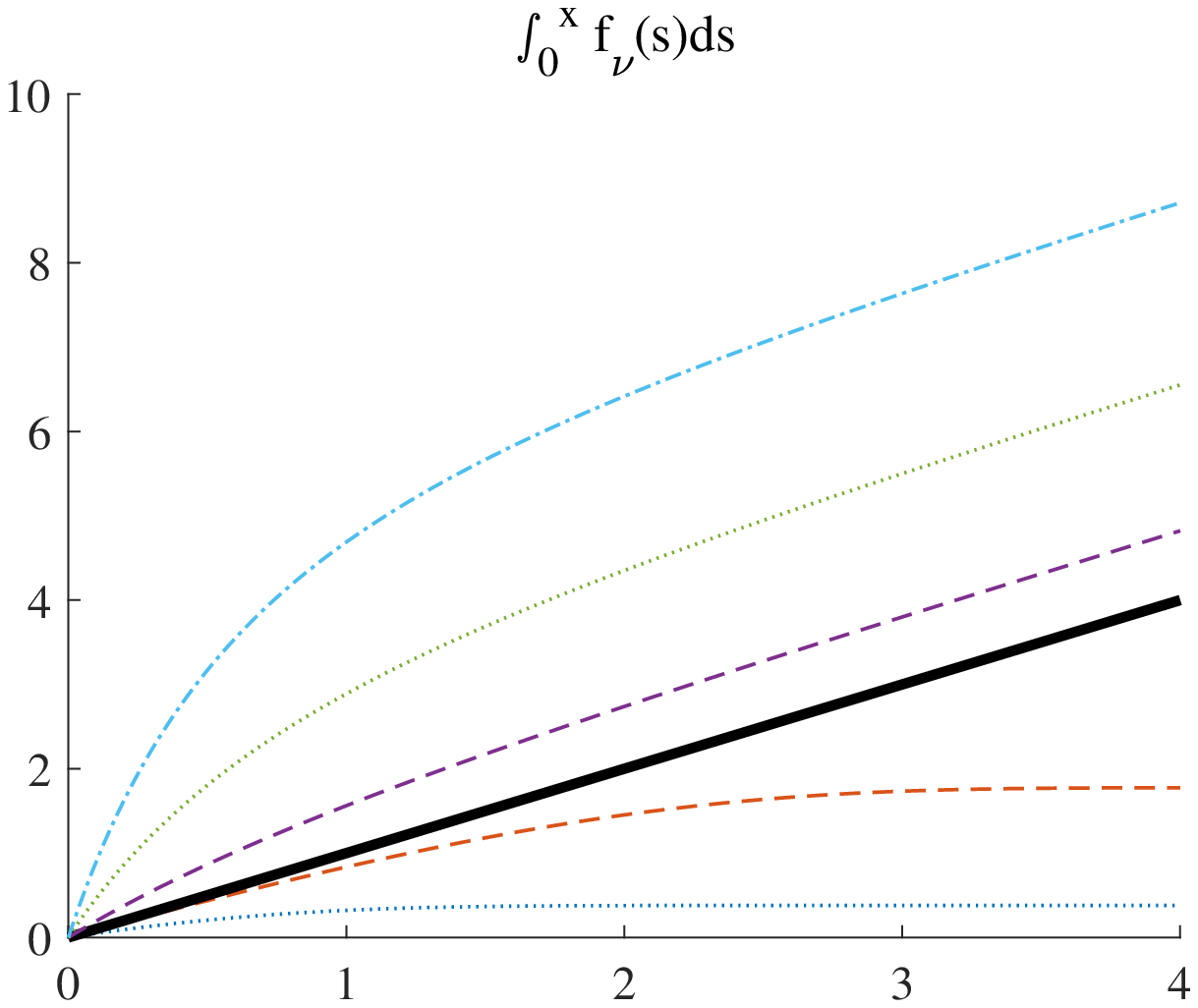}
\tiny{.}\vspace{0.1cm}
\includegraphics[height = 7cm, width = 7 cm]{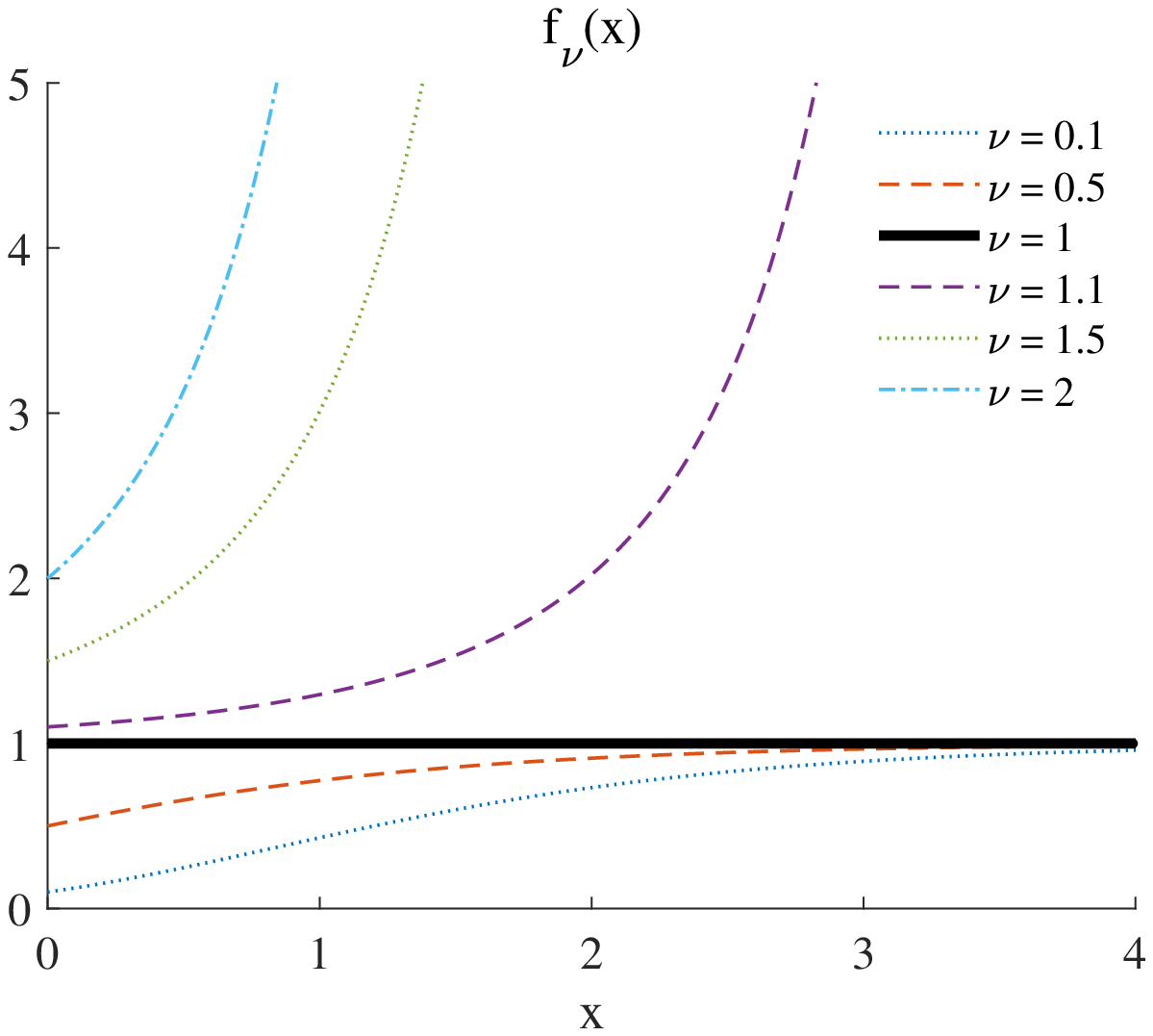}
\hspace{0.8cm}
\includegraphics[height = 7cm, width = 7 cm]{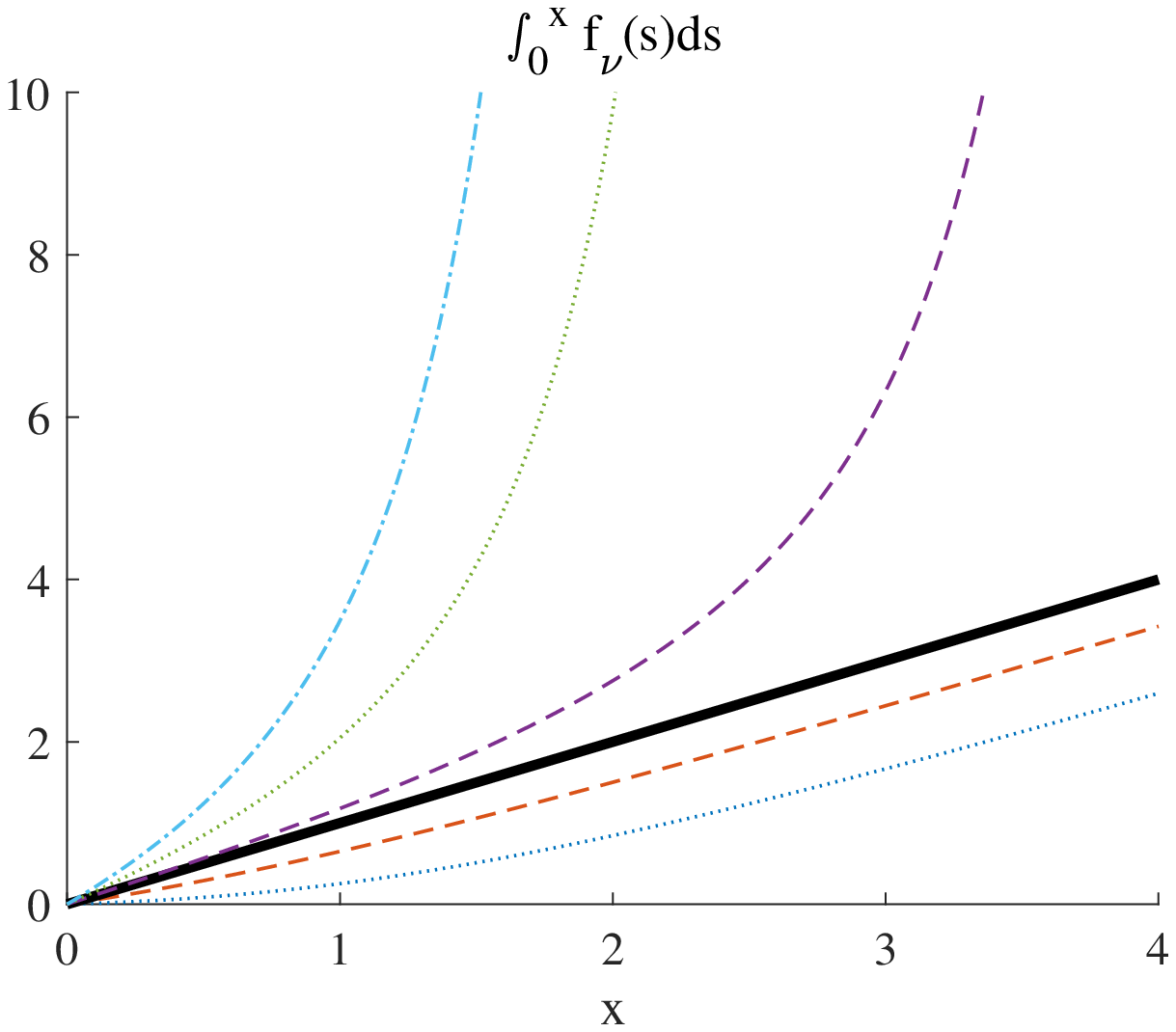}
\tiny{.}
\end{center}
\caption{Functions $f_{\nu}$ in \eqref{eq:derivative-p(x)} and $\check{u}$ in \eqref{eq:solution-variable-p} (upper row),
$f_{\nu}$ in \eqref{eq:derivative-p(x)-neg} and $\check{u}$ in \eqref{eq:solution-variable-p-neg} (lower row).
In all simulations, $A = \Lambda = 1$.}
\label{fig:variablep}
\end{figure}
%


\newpage

\section{Connections with nonlinear diffusion problems}
\setcounter{theorem}{0}
\setcounter{equation}{0}

\label{sec:diffusion}

\noindent
We follow the presentation in Lundstr\"om \cite{avhandlingen} and let $u$ denote the density of some quantity in equilibrium,
let $\Omega$ be a domain and
$E \subset \Omega$ be a $C^1$-domain so that the divergence theorem can be applied.
Due to the equilibrium, the net flux of $u$ through $\partial E$ is zero, that is
\begin{align*}
\oint\limits_{\partial E} \langle \boldsymbol F, \boldsymbol n \rangle \, ds \,= \,0,
\end{align*}

\noindent
where $\boldsymbol F$ denotes the flux density, $\boldsymbol n$ the normal to $\partial E$ and $ds$ is the surface measure.
The divergence theorem gives
\begin{align*}
\int\limits_{E} \nabla \cdot \boldsymbol F \,dx\, =\, \oint\limits_{\partial E} \langle \boldsymbol F, \boldsymbol n \rangle \,ds \,= \,0.
\end{align*}

\noindent
Since $E$ was arbitrary, we conclude
\begin{align}\label{eq:lalalalala}
\nabla \cdot \boldsymbol F\, =\, 0 \quad \mbox{in} \quad \Omega.
\end{align}

\noindent
In many situations it is physically reasonable to assume that the flux vector $\boldsymbol F$ and the gradient $\nabla u$ are related by a power-law of the form
\begin{align}\label{eq:power-law}
\boldsymbol F \,=\, - c \,|D u|^{\,q}\, D u,
\end{align}

\noindent
for some factor $c$ and exponent $q$, which may depend on space as well.
One reason is that flow is usually from regions of higher concentration to regions of lower concentration.
From this assumption, with $q = p - 2$, and from \eqref{eq:lalalalala}, we obtain the $p$-Laplace equation
\begin{align*}
\nabla \cdot (|D u|^{p-2} D u)\, = \,0 \quad \mbox{in} \quad \Omega.
\end{align*}
The linear case $p = 2$ in \eqref{eq:power-law} arises as a physical law in the following:
If $u$ denotes a chemical concentration, then it is the well known Fick's law of diffusion,
if $u$ denotes a temperature, then it is Fourier's law of heat conduction,
if $u$ denotes electrostatic potential, it is Ohm's law of electrical conduction, and
if $u$ denotes pressure, then it is Darcy's law of fluid flow through a porous media.
A problem involving the nonlinear case $p \neq 2$ is fast/slow diffusion of sandpiles, see Aronsson--Evans--Wu \cite{AEW96}.
In that case $p$ is very large and $u$ models the height of a growing sandpile.
If $|D u| > 1 + \delta$ for some $\delta > 0$, then $|D u|^{p-2}$ is very large,
and hence the transport of sand is also large,
and if $|D u| < 1 - \delta$, then $|D u|^{p-2}$ is very small.
Therefore, when adding sand particles to a sandpile,
they accumulate as long as the slope of the pile does not exceed one.
If the slope exceeds one, then the sand becomes unstable and instantly slides.
Other application in which \eqref{eq:power-law} arises with $p\neq 2$ is Hele-Shaw flow of power-law fluids (Aronsson--Janfalk \cite{AJ92}, Fabricius--Manjate--Wall \cite{FMW21}) and electro-rheological fluids (Harjulehto--H\"ast\"o--L\^e--Nuortio \cite{hhn}).
When properties of the quantity under investigation depends on space we may
model it by a variable exponent $p = p(x)$ in \eqref{eq:power-law} and thus enter
equations of type \eqref{eq:intro-p(x)} studied in Section \ref{sec:applic}.

We will now discuss the problem under investigation from the point of a diffusion problem.
Indeed, we will briefly explain, through spatially dependent diffusion,
why parts of our results presented in Theorem \ref{thm:p(x)} holds true.
Suppose that $u$ denotes the density of some quantity at equilibrium in the $n$-dimensional halfspace $\{x_n > 0\}$
and that \eqref{eq:power-law} holds with a variable exponent $p(x_n), 1 < p(x_n) < \infty$.
Assume also that $u = 0$ on the boundary $x_n = 0$, and at some $x_n = a > 0$ we assume that $u(x) > 0$.
We conclude that then $u$ satisfies the $p(x)$-Laplace equation \eqref{eq:intro-p(x)} in the halfspace and that our results applies.
We simplify by further assuming that concentration $u(x)$ is independent of $x'$-directions.
Since $|Du|$ must be positive there is a flux of $u$, independent of $x'$,
flowing perpendicular through the plane at $x_n = a$ toward the boundary $x_n = 0$.
Due to the equilibrium, the flux must be independent also of $x_n$ and is therefore constant through the halfspace. Since the problem is herefrom independent of $x'$, we drop the index and write in the following $x = x_n$.
\begin{figure}[!hbt]
\begin{center}
\includegraphics[height = 6.5cm, width = 7 cm]{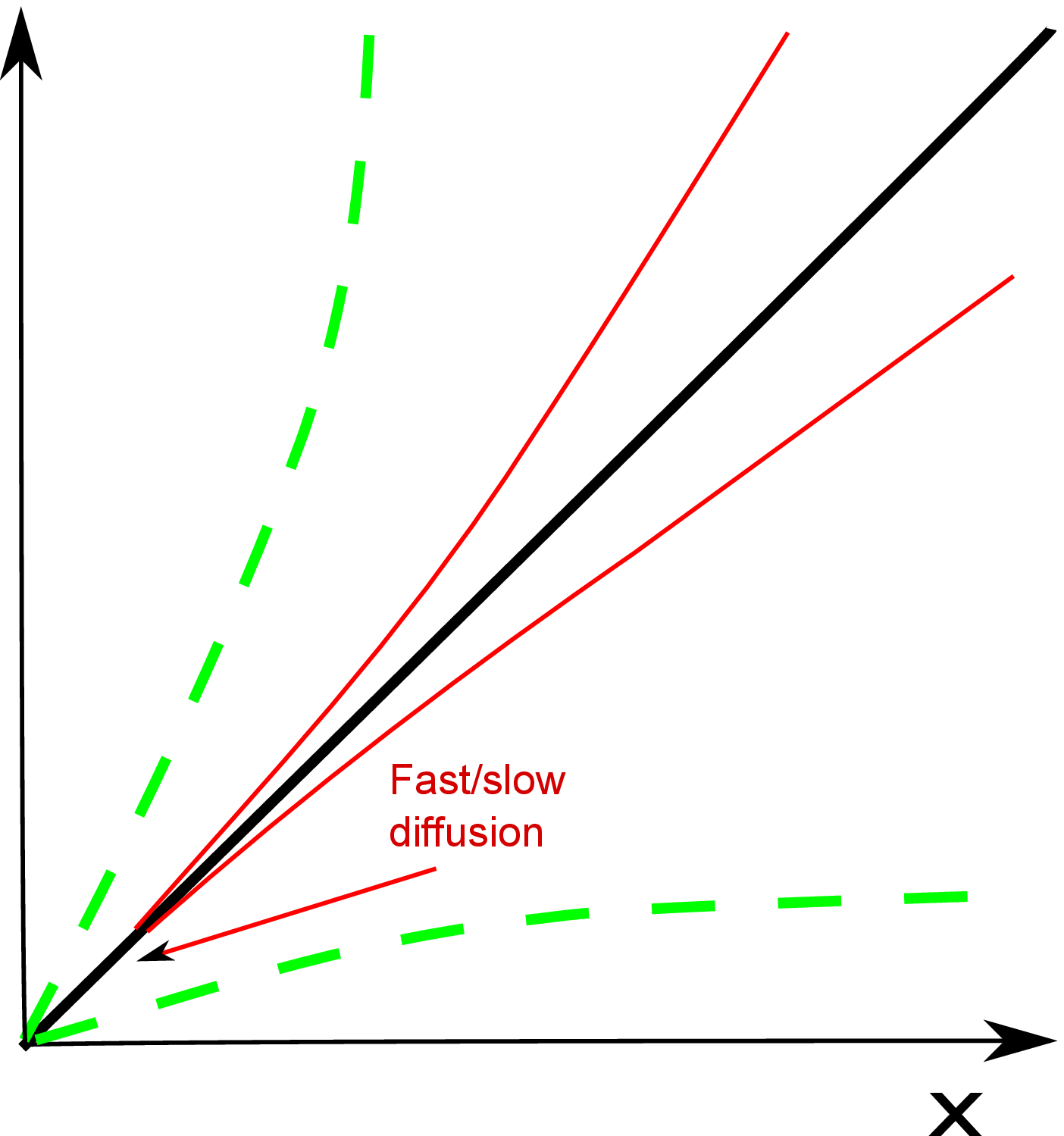}
\hspace{0.8cm}
\includegraphics[height = 6.5cm, width = 7 cm]{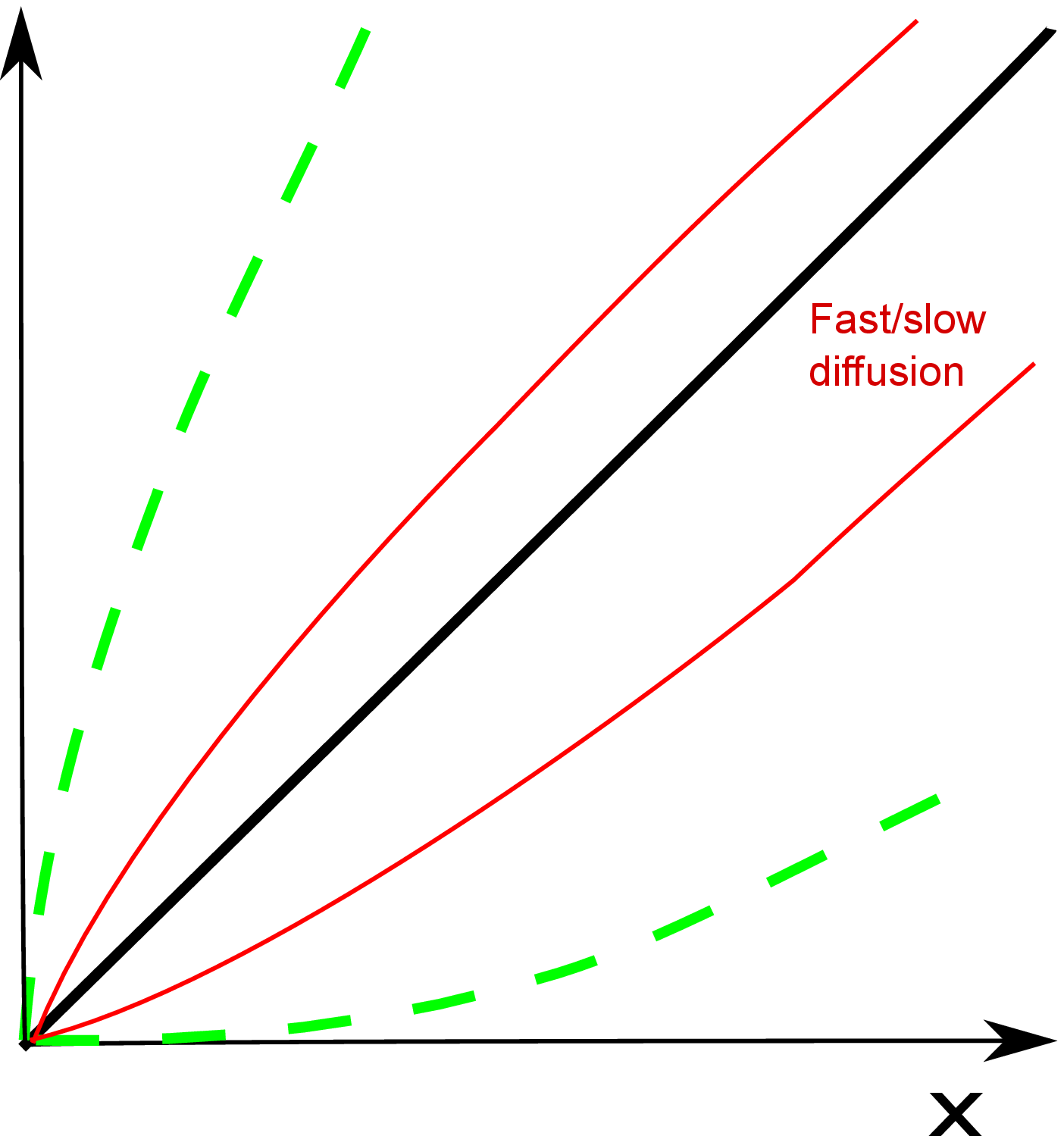}
\end{center}
\caption{Examples of how the concentration $u$ may depend on $x$ for decreasing exponents (left) and increasing exponent (right).
The slope explodes or vanish as $p(x) \to 1$: The green dashed curves correspond to an exponent $p(x)$ that approaches 1 as $x$ increases (left) and as $x \to 0$ (right). The slope approaches 1 as $p(x) \to \infty$, i.e. fast/slow diffusion: The red solid curves correspond to an exponent $p(x)$ that becomes very large as $x \to 0$ (left) and as $x$ increases (right). }
\label{fig:XXX}
\end{figure}
%

{\bf Suppose that $p(x)$ is decreasing.}
As the flux of $u$, given by assumption \eqref{eq:power-law}, is constant, the concentration $u$ must be convex (upwards) if $|D u| = u' > 1$.
Indeed, if $u' > 1$ near the boundary we locally have that \eqref{eq:power-law} yields
flux ${\bf F} = - c\, (u')^{p(x)-1}$ and since $p(x)- 1 > 0$ is decreasing it follows that $u'$ must be increasing.
A similar reasoning explains that if $u' = 1$ somewhere then the flux ${\bf F} = c$ implying
$u(x) = x$, and if $u' < 1$ then $u$ must be concave.
Figure \ref{fig:XXX}(left) shows examples of how the concentration $u(x)$ may depend on $x$ for two different decreasing exponents.
We remark that if $p(x)$ becomes very large near the boundary then $u'$ must be very close 1 there,
otherwise the flux becomes zero or infinity -- that is fast/slow diffusion (red solid curve). Similarly, if $p(x)$ comes close to 1 as we move into the domain then $u'$ must grow fast if $u'$ ever was larger than 1 along the curve in order to keep the flux constant (green dashed curve).
Finally, we realize that if $p(x)$ becomes constant then $u'$ becomes constant
(recall that $u(x) = c x$ is $p$-harmonic when $p = constant$).

{\bf Suppose now instead that $p(x)$ is decreasing.}
Reasoning as in the former case we realize that we may switch our conclusions made near the boundary in the former case with those made further away into the domain. Thus fast/slow diffusion may occur away from the boundary and in such a case the slope of $u(x)$ must approach 1.
If $p(x)$ approaches 1 near the boundary then $u'$ must explode there,
see Figure \ref{fig:XXX}(right).

Returning to \eqref{eq:solution-variable-p}, \eqref{eq:solution-variable-p-neg} and \eqref{eq:variable-p(x)-one-dimension} in Section \ref{sec:applic} we find that the
one-dimensional $p(x)$-Laplace equation yields
\begin{align*}
\Delta_{p(x)} u(x) = (p(x) - 1)u''(x) + \log |u'(x)|  p'(x) u'(x) = 0,
\end{align*}
and with the decreasing exponent 
$$
p(x) = 1 + Me^{-A x}, 
$$
where $M>0, A>0$ are a constants, the solution yields
\begin{align}\label{eq:sista1}
u(x) = \frac1A \left\{
\begin{array}{ll}
-E_i\left(\log\nu\right) + E_i\left(e^{A x}\log\nu\right)   & \text{if $\nu \neq 1$},\\
x  & \text{if $\nu = 1$}.
\end{array}
\right.
\end{align}
Similarly, with the increasing exponent
$$
p(x) = 1 + Me^{A x}
$$
the solution yields
\begin{align}\label{eq:sista2}
u(x) = \frac1A \left\{
\begin{array}{ll}
E_i\left(\log\nu\right) - E_i\left(e^{-A x}\log\nu\right)   & \text{if $\nu \neq 1$},\\
x  & \text{if $\nu = 1$}.
\end{array}
\right.
\end{align}
With $A = \lambda = \Lambda = 1$, solution curves for decreasing exponent in \eqref{eq:sista1} are plotted in Figure \ref{fig:variablep}(upper right) (below line $u = x$) and (lower right) (above line $u = x$),
and solution curves for increasing exponent in \eqref{eq:sista2} are plotted in Figure \ref{fig:variablep}(upper right) (above line $u = x$) and (lower right) (below line $u = x$).
Compare the structure of these curves to those in Figure \ref{fig:XXX} with
properties of the exponent $p(x)$ in mind. \\ 

\noindent
{\bf Acknowledgement.}
This work was partially supported by the Swedish research council grant 2018-03743.




\begin{thebibliography}{A}

\bibitem{A14}
Adamowicz T.
\emph{Phragm\'en-Lindel\"of theorems for equations with nonstandard growth},
Nonlinear Analysis: Theory, Methods and Applications {\bf 97}, (2014), 169--184.

\bibitem{A37}
Ahlfors L.,
\emph{On Phragm\'en-Lindel\"of's principle},
Trans. Amer. Math. Soc. {\bf 41}, (1937), 1--8.

\bibitem{ASS12}
Armstrong S. N., Sirakov B., Smart C. K.
\emph{Singular solutions of fully nonlinear elliptic equations and applications.}
Archive for Rational Mechanics and Analysis,  {\bf 205}(2), (2012), 345--394.

\bibitem{AEW96}
Aronsson G., Evans L. C., Wu Y.
\emph{Fast/slow diffusion and growing sandpiles.}
Journal of Differential Equations, {\bf 131} (2), (1996), 304--335.

\bibitem{AJ92}
Aronsson G., Janfalk U.
\emph{On Hele-Shaw flow of power-law fluids}.
European Journal of Applied Mathematics {\bf 3.4}, (1992) 343--366.

\bibitem{AJ17}
Avelin B., Julin V.
\emph{A Carleson type inequality for fully nonlinear elliptic equations with non-Lipschitz drift term.}
Journal of Functional Analysis {\bf 272.8}, (2017), 3176--3215.

\bibitem{Bhatt05}
Bhattacharya T.,
\emph{On the behaviour of infinity-harmonic functions on some special unbounded domains},
Pacific Journal of Mathematics {\bf 219}, no 2, (2005), 237--253.












\bibitem{BM20}
Braga J. E. M., Moreira D.
\emph{Classification of Nonnegative g-Harmonic Functions in Half-Spaces.}
Potential Analysis, (2020), 1--19.

\bibitem{CC95}
Caffarelli L.A., Cabre X.,
\emph{Fully nonlinear Elliptic equations}.
American Mathematical Society Colloquium Publications, 43.
American Mathematical Society, Providence, RI, 1995.

\bibitem{CDV07}
Capuzzo-Dolcetta I., Vitolo A.,
\emph{A qualitative Phragm\'en-Lindel\"of theorem for fully nonlinear elliptic equations},
J. Differential Equations {\bf 243}, no 2, (2007), 578--592.

\bibitem{clr}
Chen Y., Levine S., Rao M.
\emph{Variable exponent, linear growth functionals in image restoration}.
SIAM J. Appl. Math. {\bf 66} (4) (2006), 1383--1406

\bibitem{CIL92}
Crandall M. G., Ishii H., Lions P.-L.
\emph{User's guide to viscosity solutions of second order partial differential equations}, Bulletin of the American Mathematical Society,
\textbf{27} (1992), 1--67.

\bibitem{dr}
Diening L., R\r u\v zi\v cka M.
\emph{Strong solutions for generalized Newtonian fluids.}
J. Math. Fluid Mech. {\bf 7}  (2005), 413--450

\bibitem{FMW21}
Fabricius J., Manjate S., Wall P.
\emph{On pressure-driven hele-shaw flow of power-law fluids}
preprint 2021.

\bibitem{G52}
Gilbarg D.,
\emph{The Phragm\'en-Lindel\"of theorem for elliptic partial differential equations},
J. Rational Mech. Anal. {\bf 1}, (1952), 411--417.

\bibitem{GM14}
Granlund S., Marola N.
\emph{Phragm\'en-Lindel\"of theorem for infinity harmonic functions},
Commun. Pure Appl. Anal. {\bf 14} (1), (2016), 127--132.

\bibitem{hhn}
Harjulehto P., H\"ast\"o P., L\^e U. V., Nuortio M.
\emph{Overview of differential equations with non-standard growth}.
Nonlinear Analysis: Theory, Methods and Applications, {\bf 72}(12), (2010), 4551--4574.

\bibitem{H64}
Herzog J. O.
\emph{Phragmen--Lindel\"of Theorems for Second Order Quasi-Linear Elliptic Partial Differential Equations},
Proceedings of the American Mathematical Society
{\bf 15}, No. 5, (1964), 721--728.

\bibitem{H52}
Hopf E.
\emph{Remark on a preceding paper of D. Gilbarg},
J. Ration. Mech. Anal. {\bf 1}, (1952), 419--424













\bibitem{Horgan}
Horgan C. O.,
\emph{Decay estimates for boundary-value problems in linear and nonlinear continuum mechanics},
in: Mathematical Problems in Elasticity,
in: Ser. Adv. Math. Appl. Sci., {\bf 38}, World Sci. Publ, River Edge, NJ, (1996), 47--89.

\bibitem{JL03}
Jin Z., Lancaster K.,
\emph{A Phragm\'en-Lindel\"of theorem and the behavior at infinity of solutions of non-hyperbolic equations},
Pacific journal of mathematics {\bf 211}, no 1, (2003), 101--121.

\bibitem{J13}
Julin V.
\emph{Generalized Harnack inequality for nonhomogeneous elliptic equations.}
Archive for Rational Mechanics and Analysis {\bf 2.216} (2015), 673--702.

\bibitem{JJ12}
Julin V., Juutinen P.
\emph{A new proof for the equivalence of weak and viscosity solutions for the p-Laplace equation},
Communications in Partial Differential Equations {\bf 37.5} (2012), 934--946.

\bibitem{JLM01}
Juutinen P., Lindqvist P., Manfredi J. J.,
\emph{On the equivalence of viscosity solutions and weak solutions for a quasi-linear equation},
SIAM journal on mathematical analysis {\bf 33}, no 3, (2001), 699--717.




















\bibitem{JLP10}
Juutinen P., Lukkari T., Parviainen M.
\emph{Equivalence of viscosity and weak solutions for the $p(x)$-Laplacian},
Annales de l'Institut Henri Poincare (C) Non Linear Analysis. {\bf 27}. No. 6. Elsevier Masson, 2010.
1471--1487.

\bibitem{KN09}
Koike S., Nakagawa K.
\emph{Remarks on the Phragmén-Lindelöf theorem for viscosity solutions of fully nonlinear PDEs with unbounded ingredients}.
Electronic Journal of Differential Equations (EJDE)[electronic only] 2009 (2009): Paper-No.

\bibitem{K93}
Kurta V. V.,
\emph{Phragm\'en-Lindel\"of theorems for second-order quasilinear elliptic equations},
(Russian)
Ukrain. Mat. Zh. {\bf 44}, no 10 (1992), 1376--1381;
translation in Ukrainian Math. J. {\bf 44}, no 10 (1992),
1262--1268 (1993).

\bibitem{Qnew}
M. Leseduarte, M. Carme, R Quintanilla,
\emph{Phragmén-Lindelöf alternative for the Laplace equation with dynamic boundary conditions}
Journal of applied analysis and computation 7.4 (2017): 1323--1335.


\bibitem{L85}
Lindqvist P.,
\emph{On the growth of the solutions of the differential equation $\nabla \cdot (|\nabla u|^{p-2}\nabla u) = 0$ in $n$-dimensional space},
Journal of Differential Equations, {\bf 58}, (1985), 307--317.

\bibitem{LW15}
Lundberg E., Weitsman A.
\emph{On the growth of solutions to the minimal surface equation over domains containing a halfplane.}
Calculus of Variations and Partial Differential Equations, {\bf 54}(4), (2015), 3385--3395.

\bibitem{avhandlingen}
Lundstr\"om N. L. P.,
{\em $p$-harmonic functions near the boundary,}
Doctoral Thesis, ISSN 1102-8300, ISBN 978-91-7459-287-0, Ume\aa $\;$ 2011.

\bibitem{L16}
Lundstr\"om N. L. P.,
\emph{Phragmén-Lindelöf Theorems and p-harmonic Measures for Sets Near Low-dimensional Hyperplanes},
Potential Analysis, {\bf 44}, (2016), 313--330.

\bibitem{LOT20}
Lundstr\"om N. L. P., Olofsson M., Toivanen O.
\emph{Strong maximum principle and boundary estimates for nonhomogeneous elliptic equations.}
arXiv preprint arXiv:2005.03338 (2020).

\bibitem{LS21}
Lundstr\"om N. L. P., Singh J.
\emph{Estimates of p-harmonic functions in planar sectors.}
arXiv preprint arXiv:2111.02721 (2021).

\bibitem{MO17}
Medina M., Ochoa P.
\emph{On viscosity and weak solutions for non-homogeneous p-Laplace equations}
Advances in Nonlinear Analysis {\bf 8} (2017), 468--481.

\bibitem{M71}
Miller K.
\emph{Extremal barriers on cones with Phragmen--Lindel\"of theorems and other applications}. Annali di Matematica Pura ed Applicata, 90(1), (1971), 297--329.

\bibitem{PL08}
Phragm\'en E., Lindel\"of E.,
{\em Sur une extension d'un principe classique de l'analyse et sur quelques propri\'et\'es des functions monog\'enes dans le voisinage d'un point singulier,}
Acta Math. {\bf 31}, no 1, (1908), 381--406.

\bibitem{Q93}
Quintanilla R.,
\emph{Some theorems of Phragm\'en-Lindel\"of type for nonlinear partial differential equations},
Publ. Mat {\bf 37}, (1993), 443--463.

\bibitem{S54}
Serrin J.,
\emph{On the Phragm\'en-Lindel\"of principle for elliptic differential equations},
J. Rational Mech. Anal. {\bf 3}, (1954), 395--413.

\bibitem{V04}
Vitolo A.,
\emph{On the Phragm\'en-Lindel\"of principle for second-order elliptic equations},
J. Math. Anal. Appl. {\bf 300}, no 1, (2004), 244--259.





\end{thebibliography}
\end{document}